\documentclass[ejs]{imsart}

\RequirePackage{amsthm,amsmath}
\RequirePackage[numbers]{natbib}
\RequirePackage[colorlinks,citecolor=blue,urlcolor=blue]{hyperref}

\usepackage{amsmath,amsfonts,amssymb,amsthm,booktabs,color,epsfig,graphicx,hyperref,url}

\usepackage{dsfont} 
\usepackage{chngcntr}
\usepackage{apptools}
\usepackage[ruled,vlined,linesnumbered]{algorithm2e}
\usepackage{relsize}
\usepackage{array,multirow,makecell}
\usepackage{hhline}
\usepackage{bm}
\usepackage{placeins}
\usepackage{rotating}
\usepackage{caption,subcaption}
\usepackage{adjustbox}

\doi{10.1214/154957804100000000}
\pubyear{0000}
\volume{0}
\firstpage{1}
\lastpage{1}

\startlocaldefs

\theoremstyle{plain}
\newtheorem{theorem}{Theorem}

\newtheorem{proposition}{Proposition}

\newtheorem{corollary}{Corollary}

\theoremstyle{definition}

\newtheorem{remark}{Remark}

\def \m {\mathbf}
\def \R {\mathbb R}
\def \beg{\begin{eqnarray}}
\def \en{\end{eqnarray}}
\def \be*{\begin{eqnarray*}}
\def \e*{\end{eqnarray*}}
\def \di{\displaystyle}
\def \bit{\begin{itemize}}
\def \eit{\end{itemize}}
\def \E{\mathbb E}
\def \N{\mathbb N}
\def \I{\mathbb I}
\def \cqfd{\blacksquare}
\def \si{\sigma}

\def \V {\mathbb{V}}

\renewcommand{\P}{\mathbb{P}}
\makeatletter
\newcommand{\RemoveAlgoNumber}{\renewcommand{\fnum@algocf}{\AlCapSty{\AlCapFnt\algorithmcfname}}}
\newcommand{\RevertAlgoNumber}{\algocf@resetfnum}
\makeatother
\def \w{\widehat}
\def \t{\tilde}
\def \eps{\epsilon}
 \def \cqfd{\hspace*{10.6cm}$\blacksquare$}

\def\argmax{\mathop{\mbox{\sl\em argmax}}}

\endlocaldefs

\begin{document}

\begin{frontmatter}

\title{Smooth test for equality of copulas}
\runtitle{Smooth test for equality of copulas}

\begin{aug}

\author[1]{\fnms{Yves Ismaël} \snm{Ngounou Bakam}\thanksref{t1}\corref{}\ead[label=e1]{yves.ngounou@ensai.fr}}
\address[1]{Univ Rennes, Ensai, CNRS, CREST - UMR 9194, F-35000 Rennes, France \\\printead{e1} }
\and
\author[1,2]{\fnms{Denys} \snm{Pommeret}\ead[label= e2]{denys.pommeret@univ-amu.fr}}
\address[2]{ISFA, Univ Lyon, UCBL, LSAF EA2429, F-69007, Lyon, France\\\printead{e2}}
\thankstext{t1}{corresponding author}
\runauthor{Ngounou Bakam $\&$ Pommeret}
\end{aug}

\begin{abstract}

A smooth test to simultaneously compare $K$ copulas, where $K \geq 2$ is proposed. The $K$ observed populations can be paired, and the test statistic is constructed based on the differences between moment sequences, called copula coefficients. These coefficients characterize the copulas, even when the copula densities may not exist.
The procedure employs a two-step data-driven procedure. In the initial step, the most significantly different coefficients are selected for all pairs of populations. The subsequent step utilizes these coefficients to identify populations that exhibit significant differences.
To demonstrate the effectiveness of the method, we provide illustrations through numerical studies and application to two real datasets.
%

\end{abstract}




\tableofcontents

\end{frontmatter}


\section{Introduction and motivations}
Copulas have been extensively studied in the statistical literature and their field of application covers a very wide variety of areas (see for instance the book of  \cite{joe} and references therein). 
The problem of goodness-of-fit for  copulas is, therefore, an important topic and can deserve many situations as in insurance to compare the dependence between portfolios (see for instance \cite{Shi16}),  in finance to compare the dependence between indices (see for instance the book of \cite{cherubini2004copula}), in biology to compare the  dependence between genes (see \cite{kim2008copula}), in medicine to compare  diagnosis  (see for instance \cite{hoyer}), 
  or more recently in ecology to compare dependence between species (see \cite{ghosh}).

In the one-sample case, many testing methods have been proposed in the context  of parametric families of copulas
(see for instance  the review paper of \cite{genest2009}, 
or more recently \cite{Omelka}, 
\cite{can2015}, and \cite{can2020}).

 In the two-sample case,  an important reference is the work of \cite{remillard}. 
    They proposed a nonparametric test based on the integrated square difference between the empirical copulas. Their approach requires the continuity of partial derivatives of copulas which allows to obtain an approximation of the distribution under the null.  Their test is convergent and is adapted to independent as well as paired populations, and an R package {'TwoCop'} is available (see \cite{twocop}). 

When $K>2$,
\cite{quessy16} proposed an innovative  work to compare  $K$ copulas. More recently \cite{quessy21} developed a second test statistic with a very original idea based on a generalized Szekely–Rizzo inequality.
These tests are consistent and can also be used to test radial symmetry and exchangeability of copulas. However,
\cite{quessy16, quessy21}  restricted his study to the case of samples of the same size. More precisely both procedures consist
 of dividing the sample into sub-samples and testing the equality of the associated sub-copulas.
Therefore, testing the equality of copulas from independent samples cannot be achieved by these works.
Furthermore, in both cases the null distribution is intractable and the author needs a multiplier bootstrap method to implement these tests.
 Such bootstrap approach for copulas was initiated in \cite{scaillet}. 
Another  extension of \cite{remillard} is proposed in
\cite{bouzebda} when the $K$ populations are observed independently, but 
the proposed test statistic seems to work only for testing the simultaneous independence of the $K$ populations.

Recently, \cite{Ngounou} studied a nonparametric copula estimator which showed very good numerical results.
In this paper, we propose to address the problem of $K$-copulas comparison with a new approach based on such estimators.
 We do not directly compare the empirical copulas, but we compare their projections on the basis of Legendre polynomials. We restrict our study to continuous variables whose populations can be paired. 
Then it makes possible to simultaneously compare the dependence structures of various populations, such as various portfolios in insurance, as well as to compare the same population followed over several periods, such as medical cohorts. 
Moreover, the procedure is valid 
for the case of several independent samples with different sample sizes, which is important for applications and a novelty compared to the works cited above, even if the works of \cite{quessy16, quessy21} could certainly be generalized in this direction.

Our method is a  data-driven procedure derived from the Neyman's smooth tests theory (see \cite{neyman1937smooth}). These smooth tests are omnibus tests and detect any departure from the null.
In our case, we consider the orthogonal projections of the copula densities on the basis of Legendre polynomials and we compare their coefficients.
For each pair of populations, a penalized rule is introduced to select automatically the coefficients that are the most significantly different. A second penalized rule selects the number of populations to be compared. Thus the procedure is a data-driven method with two selection steps.
Under the null, due to the penalties, the rules select only one pair of populations and only one coefficient, leading to a chi-square asymptotic null distribution. Then
the test is very simple and easy to implement.
This is another major difference from the work of \cite{quessy16, quessy21} where the null distribution does not have an explicit form in general and where a multiplier bootstrap is used to calculate the p-values.
We also prove that the test procedure detects any fixed alternative and gives us information on the reject decision. More precisely, the second  penalized rule
 is calibrated to detect 
 the populations that differ most
 significantly. Then in case of rejection, we can find the pairs of populations that contributed the most to the value of the test statistic. We can also proceed to a two-by-two test to search similar populations.  In practice, we have developed an R package 'Kcop' which is available on the Comprehensive R Archive Network (CRAN)  to implement the $K$-sample procedure.

A numerical study shows the good behaviour of the test.
We apply this approach on two datasets related to biology and insurance. The first one is the very well-known Iris dataset. While this dataset is very famous there was no simultaneous comparison between the 4-dimensional dependence structures of the three species involved. We therefore propose to apply the smooth test to compare the dependence between sepals and petals, thus providing a new analysis. 
The second dataset is a large medical insurance database with possibly paired data and concerns claims from three years: 1997, 1998 and 1999. We apply the smooth test on several variables from this dataset illustrating the idea of risk pooling and price segmentation.  
All these results can be reproduced using the   'Kcop' package.

The paper is organized as follows: in Section 2 we specify the null hypothesis considered in this paper and we set up the notation.
Section 3 presents the method in the two-sample case. 
In Section 4 we extend the result to the $K$ ($K>2$) sample case and in Section 5 we proceed with the study of the convergence of the test under 
alternatives.
Section 6 is devoted to the numerical study and Section 7 contains real-life illustrations.
Section 8 discusses extensions and connections.

All proofs are located in Appendix \ref{appendproofs}. The adaptation to the dependent case is straightforward and is summarized in Appendix \ref{appendIndep}, where all results are rewritten in this context. A method for automating test parameters is available in Appendix \ref{annexTuning}. Additionally, Appendices \ref{appendlegendre} to \ref{appendempirical} contain supplementary materials, including various complements, additional simulations, and comparisons.

\section{Notation and null hypotheses}

Let $\mathbf X=(X_1,\ldots, X_p)$  be a $p$-dimensional  continuous random vector with  joint cumulative  distribution  function (cdf)  $F_{\mathbf X}$, and with unique copula defined by  
\be*
\label{cop}
C(F_{1}(x_1),\ldots, F_{p}(x_p) )  & = &
F_{\m X}(x_1,\ldots, x_p),
\e*
where $F_{j}$  denotes the marginal cdf of $X_j$. 
Writing
\be*
U_j : = F_{j}(X_j),&&{\rm  for \ }j=1,\ldots, p,
\e*
we have  for all $u_j \in [0,1]$
\be*
C(u_1,\ldots,u_p) & = & F_{\m U}(u_1,\ldots,u_p),
\e*
with $\m U=(U_1,\ldots, U_p)$. 
The copula density (if it exists)
defined by 
\be*
c(u_1,\ldots, u_p) & : = & 
\di\frac{\partial^p C( u_1,\ldots, u_p)}{\partial u_1,\ldots, \partial u_p},
\e*
coincides with the probability density function (pdf) $f_{\m U}$ 
of the vector $\mathbf U$.
Write ${\cal L}=\{L_n ; n\in\N\}$ the set of orthogonal Legendre polynomials with first terms $L_0=1$ and $L_1(x)=\sqrt{3}(2x-1)$, such
that  $L_n$ is of degree $n$ and satisfies (see Appendix \ref{appendlegendre}  for more details):
\be*
\di\int_0^1 L_j(u) L_k(u) du &  = & \delta_{jk},
\e*
where $\delta_{jk}=1$ if $j=k$ and $0$ otherwise.
The random variables $U_i$ are uniformly distributed and we have the following 
decomposition
\beg
\label{density2}
c(u_1,\ldots  ,u_p)
& = &
\di\sum_{j_1, \ldots, j_p \in \N}
\rho_{j_1,\ldots, j_p} L_{j_1}(u_1)\ldots L_{j_p}(u_p),
\en
where
\be*
\rho_{j_1,\ldots ,j_p} & = &
\E(L_{j_1}(U_1)\ldots L_{j_p}(U_p)),
\e*
as soon as $f_{\m U}$ exists and  belongs to
the space of all square-integrable
functions with respect to the  Lebesgue measure on $[0,1]^{ p}$, that is, if
\beg
\label{condition}
\di \int_0^1 \ldots \di  \int_0^1 c(u_1,\ldots, u_p)^2 du_1\ldots du_p & < & \infty .
\en
Write $\m j=(j_1,\ldots, j_p)$ and $\m 0=(0,\ldots, 0)$. We can observe that $\rho_{\m 0} =1$. Moreover,  since by
orthogonality we have $\E(L_{j_i}(U_i))=0$, for all $i=1,\ldots, p$, we see that   $\rho_{\m j}=0 $ if only one element of $\m j$ is non null. 
When the copula density exists and is square integrable, we deduce from  
(\ref{density2})  
that, for all $u_1,\ldots, u_p \in [0,1]$, 
\beg
c(u_1,\ldots, u_p)
& = &
1+ \di\sum_{\m j \in \N^p_*}
\rho_{\m j} L_{j_1}(u_1)\ldots  L_{j_p}(u_p),
\nonumber
\\
\label{colula}
C(u_1,\ldots, u_p) & = & \di u_1 \; u_2\ldots u_p+\sum_{\m j \in\N^p_*}
\rho_{\m j} I_{j_1}(u_1) \ldots I_{j_p}(u_p),
\en
where
$
I_j(u) = \di\int_{0}^{u}L_j(x) dx$,
and $\N^p_*$ stands for the set $\{\m j =(j_1,\ldots, j_p)  \in \N^p ; \m j \neq \m 0\}$.  The sequence $(\rho_{\m j})_{\m j  \in \N^p_*}$  will be referred to as the {\it copula coefficients} (as in  \cite{Ngounou}). Since $\m U$ is bounded, all copula coefficients exist.  
The following result, due to \cite{moment2} or \cite{kleiber2013multivariate}, shows that 
such a  sequence  
characterizes the copula. Moreover, it shows that assumption (\ref{condition}) is unnecessary.  
\begin{proposition}
\label{prop00}
Let $(\rho_{\m j})_{\m j  \in \N^p}$ and 
$(\rho'_{\m j})_{\m j  \in \N^p}$ be two sequence of copula coefficients associated to copulas $C$ and $C'$, respectively.  Then 
\begin{align*}
 \rho_{\m j}=\rho'_{\m j}, \ \forall \m j  \in \N^p & \Longleftrightarrow C=C'.
 \end{align*}
 \end{proposition}
Thereby, the copula is determined by its sequence of copula coefficients, a property that holds even when condition  (\ref{condition}) is not satisfied, and the copula density may not exist.
Consequently, for any continuous random vectors,  the comparison of their copulas coincides with the comparison of 
their copula coefficients. This equivalence holds true even when the random vectors lack a density or possess densities that are not square-integrable.  We will use this characterization to construct the test statistic. 

We consider $K$ continuous random vectors, namely
$$\mathbf X^{(1)} = (X^{(1)}_{1},\ldots, X^{(1)}_{p}), \ldots, \mathbf X^{(K)} = (X^{(K)}_{1},\ldots, X^{(K)}_{p}),$$
 with joint cdf  $\mathbf F^{(1)},\ldots, \mathbf F^{(K)}$, and with associated copulas  $C_1, \ldots, C_K$, respectively.
Assume that we observe $K$  iid samples from $\m X^{(1)}, \ldots, \m X^{(K)}$, possibly  paired,  denoted by
$$
(X^{(1)}_{i,1},\ldots, X^{(1)}_{i,p})_{i=1,\ldots, n_1}, \ldots,
( X^{(K)}_{i,1},\ldots, X^{(K)}_{i,p})_{i=1,\ldots, n_K}.
$$
The following assumption will be needed throughout the paper: we assume that for all $1\leq \ell < m   \leq K$, $\min(n_{\ell}, n_m) \rightarrow \infty$, and
\beg
\label{ratiorate}
 n_{\ell}/(n_{\ell}+n_{m}) \rightarrow a_{ \ell, m},  {\rm \ with \ } 0< a_{\ell, m} < \infty.
\en
Write $\m n=(n_1,\ldots,n_K)$. Hence, it will cause no confusion if we write $\m n \rightarrow +\infty$ when all $n_i \rightarrow +\infty$, and for a series of  univariate  random variable $(Q_n)_{n\in\N}$ the notation 
$Q_n = o_{\P}(\m n)$ means that 
$Q_n = o_{\P}(n_i)$, for all $i=1,\ldots, K$.  

We  consider the problem of testing the equality
\beg
\label{equality}
H_0: && C_1 = \cdots =  C_K,
\en
 against $
H_1:$
 there  exist  $1\leq k \neq k' \leq K$  such  that   $C_k \neq  C_{k'}$. 
From
Proposition \ref{prop00}, testing the equality (\ref{equality}) remains to test the equality of all copula  coefficients, that is
\beg
\label{hypothesis2}
 H_0: \rho^{(1)}_{\m j} = & \cdots &
= \rho^{(K)}_{\m j} , \ \ \forall \m j \in \N^p_* ,
\en
 against $
 H_1:$
 there  exist  $1\leq k \neq k' \leq K$  and $\m j \neq \m j'$ such  that   $\rho^{(k)}_{\m j} \neq \rho^{(k')}_{\m j'}$,
where $\rho^{(k)}$ stands for the copula coefficients associated to $C_{k}$.

We will denote by $F^{(\ell)}_{ j}$ the marginal cdf  of the  $j$th component of $\m X^{(\ell)}$ and we write
\begin{align*}
   U^{(\ell)}_{i,j} & =
   F^{(\ell)}_{ j}(X^{(\ell)}_{i,j}).
\end{align*} 
For testing (\ref{hypothesis2}), we  estimate the copula coefficients by
\be*
 \w \rho^{(\ell)}_{j_1\ldots j_p} & = & \di\frac{1}{n_{\ell}} \di\sum_{i=1}^{n_{\ell}} L_{j_1}(\w U^{(\ell)}_{i,1})\ldots L_{j_p}(\w  U^{(\ell)}_{i,p}),
\e*
{\rm where }
 $\w  U^{(\ell)}_{i,j}  =
   \w F^{(\ell)}_{ j}(X^{(\ell)}_{i,j})$,
and $\w F$ denotes the empirical distribution function associated to $F$. Such estimators $ \w \rho^{(\ell)}_{j_1\ldots j_p}$ have been extensively studied in \cite{Ngounou} where it is shown their excellent behavior. 
Considering the null hypothesis $ H_0$ as expressed in  (\ref{hypothesis2}), our  test procedure is based on the   sequences
of  differences
\be*
r^{(\ell,m )}_{\m j} :=  \w \rho^{(\ell)}_{\m j}-\w { \rho}^{(m)}_{\m j},  &&{\rm for \ } 1\leq  \ell \leq m \leq K, {\rm \ and \ } \m j \in \N^p_*,
\e*
with the convention that  $r^{( \ell, m)}_{\m j}=0$ when only one component of $\m j$ is different from zero. This is due to the orthogonality of the Legendre polynomials, leading 
$\rho^{(\ell)}_{\m j}=\rho^{(m)}_{\m j}=0$ in such cases. 

In order to select automatically the  number of copula coefficients, 
for any vector $\m j=(j_1,\ldots, j_p)$, we will denote by
$$\|\m j\|_1=|j_1|+\cdots + |j_p|,$$
the $L^1$ norm and for any integer $d>1$,  we write
\be*
{\cal S}(d) = \{ \m j\in\N^p;\, \|\m j\|_1=d \text{ and }\exists \,k\neq k' \text{ such that }
j_k > 0\text{ and }j_{k'} > 0\}.
\e*
The set ${\cal S}(d)$ contains all non null positive  integers $\m j=(j_1,\ldots,j_p)$ with $L^1$ norm equal to $d$  and such that $j_k <d$, for  all $k=1,\ldots, p$.
We will denote by  $ c(d):=\binom{d+p-1}{d}-p $ the cardinality of ${\cal S}(d)$
and we introduce a lexicographic order   
on $\m j\in{\cal S}(d)$  as follows:
\be*
\m j=(d-1,1,0,\ldots, 0)  & \Rightarrow & ord_{{}}(\m j,d) = 1
\\
\m j=(d-1,0,1,\ldots, 0)   & \Rightarrow & ord_{{}}(\m j,d) = 2
\\
&\ldots &
\\
\m j=(0,\ldots, 0,2,d-2) & \Rightarrow & ord_{{}}(\m j,d)= c(d)-1
\\
\m j=(0,\ldots, 0,1,d-1)  & \Rightarrow & ord_{{}}(\m j,d)= c(d).
\e*
This order will be used to compare successively the copula coefficients. 

\section{Two-sample case}
We first consider the two-sample case when $K=2$ to detail the construction of the test statistics. We want
to test
\be*
 H_0: \rho^{(1)}_{\m j}
= \rho^{(2)}_{\m j} , \ \ \forall \m j \in \N^p_*.
\e*
We  restrict our attention to the 
iid case, 
the paired case with $n_1 = n_2$ being briefly described  in  Appendix \ref{appendIndep}.
To compare the copulas associated with $\m X^{(1)}$ and $\m X^{(2)}$, we introduce  a series of statistics derived from the differences between their copula  coefficients. Specifically,   
for $1\leq k \leq  c(2)$,  we define
\begin{eqnarray}\label{Tkindep}
{T}^{(1, 2)}_{2,k}&:=& \di\frac{n_{1}  n_2}{n_{1}+ n_2} \di\sum_{\m j\in{\cal S}(2); ord_{{}}(\m j,2)\leq k}   (r^{(1, 2)}_{\m j})^2,
\end{eqnarray}
and, for $d>2$ and $1\leq k \leq c(d)$,
\begin{eqnarray}
\label{Tk2indep}
{T}^{(1, 2)}_{d,k}&:=& T^{(1, 2)}_{d-1,c(d-1)} + \di\frac{n_{1}  n_2}{n_{1}+ n_2} \di\sum_{\m j\in {\cal S}(d); ord_{{}}(\m j,d) \leq k }  ( r^{(1, 2)}_{\m j})^2.
\end{eqnarray}
These statistics are embedded and  we have for $2\leq k < c(d)$,
\be*
T^{(1, 2)}_{d,k}
& = &
 \di\frac{n_{1}  n_2}{n_{1}+ n_2} \left(
 \di\sum_{u=2}^{d-1}
 \di\sum_{\m j \in{\cal S}(u)}(r^{(1, 2)}_{\m j})^2+
  \di\sum_{\m j\in {\cal S}(d); ord_{{}}(\m j,d) \leq k }  ( r^{(1, 2)}_{\m j})^2\right).
\e*
It follows that
\be*
T^{(1, 2)}_{2,1} \leq T^{(1, 2)}_{2,2} \leq \cdots \leq T^{(1, 2)}_{2,c(2)} \leq  T^{(1, 2)}_{3,1} \leq  \cdots \leq T^{(1, 2)}_{d,c(d)} \leq T^{(1, 2)}_{d+1,1} \leq \cdots .
\e*
Each statistic $T^{(1, 2)}_{d,k}$ contains information enabling the comparison  of the copula   coefficients $\rho_{\m j}^{(1)}$ and $\rho_{\m j}^{(2)}$ up to the norm
$\|\m j\|_1 = d$ and $ ord_{{}}(\m j,d)=k$. Consequently, for a large value of $d$, it will be possible to compare the coefficient of high orders using  $r_{\m j}^{(1,2)}$, while the parameter $k$ allows the exploration of all values of $\m j$ for the given order.
To simplify notation, we write such a sequence of statistics as
\be*
V^{(1, 2)}_1 = T^{(1, 2)}_{2,1} ; \ V^{(1, 2)}_2 =  T^{(1, 2)}_{2,2} ; \ldots \ V^{(1, 2)}_{c(2)} =  T^{(1, 2)}_{2,c(2)} ; \  V^{(1, 2)}_{c(2)+1} = T^{(1, 2)}_{3,1} \ldots \
\e*
By construction, for all integer $k>0$, each statistic $V^{(1, 2)}_k $ is a sum of $k$ elements. More precisely  there exists a set $ {\cal H}(k) \subset \N_*^p$,  with $card({\cal H}(k))=k$,  such that
\begin{align}
    \label{sumV}
V^{(1, 2)}_k &=
 \di\frac{n_{1}  n_2}{n_{1}+ n_2} \di \sum_{\m j \in{\cal  H}(k)}
  (r^{(1, 2)}_{\m j})^2.
\end{align}
It can be observed that if $\m j $ belongs to ${\cal H}(k)$ then $\|\m j\|_1 \leq k$.
Moreover, we have the following relation: for all $k\geq 1$ and $j=1,\ldots, c(k+1)$
\begin{align*}
V^{(1, 2)}_{c(1)+c(2)+\cdots+c(k)+j}
& =
T^{(1, 2)}_{k+1,j},
\ \ \ \ {\rm with \ the \ convention \ } c(1)=0.
\end{align*}

Notice that we need to compare all copula coefficients and then let $k$ tend to infinity to detect all possible alternatives. 
However, choosing a too large value for $k$ can lead to a dilution of the test's power. 
Following  \cite{KL95}, we suggest a  data-driven procedure  to automatically select  the number of coefficients
to test the hypothesis $H_0$.
For this purpose, we set
\begin{align} \label{rule}
D(\m n) & :=  \min\big\{\argmax_{1\leq k \leq d(\m n)} (  V^{(1, 2)}_{k} - k p_{\m n})  \big\},
\end{align}
where
$p_{\m n}$  and $d(\m n)$ tend to $ +\infty$, as $n_{1}, n_2 \to +\infty$, $k p_{\m n}$ being a penalty term which penalizes the embedded statistics proportionally to the number of copula coefficients used.
Roughly speaking, $D(\m n)$   automatically selects the coefficients that exhibit the most significant differences.

Therefore, the data-driven test statistic that we use to
compare $C_1$ and $C_2$ is $V^{(1, 2)}_{D(\m n)}$. We consider the following rate for penalty term: 

\begin{description}
\item{{\bf (A)}}
$d(n_i)^{(p+5)} = o(p_{\m n})$, {\rm \ for \ } i=1,2.
\end{description}

Our first result shows that under the null the least penalized statistic will be selected, that is, the first one.
\begin{theorem}
\label{thm1_ind}
If  {\bf (A)} holds, then, under $ H_0$, ${D(\m n)}$ converges in probability towards 1 as $n_{1},n_2 \rightarrow +\infty$.
\end{theorem}

It is worth noting that under the null, the asymptotic distribution of the statistic $V^{({1},2)}_{D(n)}$ coincides with the asymptotic distribution of $V_{1}^{(1, 2)}=T^{(1, 2)}_{2,1}= \di\frac{n_{1}  n_2}{n_{1}+ n_2}(r_{\m j}^{(1,2)})^2$, with $\m j = (1,1,0,\ldots,0)$.  In that case, we simply have
\begin{align*}
    r^{(1, 2)}_{\m j}
    & =
    \di\frac{1}{n_1}
    \di\sum_{i=1}^{n_1}
L_1(\w U^{(1)}_{i,1}) L_1(\w U^{(1)}_{i,2})
    -
 \di\frac{1}{n_2}
    \di\sum_{i=1}^{n_2}   L_1(\w U^{(2)}_{i,1}) L_1(\w U^{(2)}_{i,2}).
\end{align*}
It follows that $T_{2,1}^{(1, 2)}$ measures  the discrepancy between
$\E(L_1( U_{1}^{(1)})L_1( U_{2}^{(1)}))$ and $\E(L_1( U_{1}^{(2)})L_1( U_{2}^{(2)}))$. 
This simply means that all other copula coefficients are not significant under the null and are therefore not selected.
Asymptotically, the null distribution reduces to that of $V_{1}^{(1, 2)}$ and is given below.
\begin{theorem}\label{thm2_ind}
Let $\textbf{j}=(1,1,0\ldots,0)$. 
Then under $H_0$, 
\be*
\di V_1^{(1, 2)}/\sigma^2(1,2) &\overset{D}{\longrightarrow}& \chi^2_1,
\e*
with 
$\sigma^2(1,2) =(1-a_{1 ,2}) \sigma^2(1)+a_{1, 2} \sigma^2(2)$,  where $a_{1,2}$ is defined in (\ref{ratiorate}), and where, for $s=1,2$,  
\be*
&&\sigma^2(s) =\V
 \Bigg(
 L_1(U^{(s)}_{1})L_1(U^{(s)}_{2})\\
&&\qquad\qquad\quad+ 2\sqrt{3}\di\int\di\int \big(\mathds{1}(X^{(s)}_{1}\leq x )-F_{1}^{(s)}(x)\big)L_{1}( F^{(s)}_{2}(y))
d F^{(s)}(x,y)
\\
&&\qquad\qquad\quad+2\sqrt{3}\di\int\di\int \big(\mathds{1}(X^{(s)}_{2}\leq y )-F_{2}^{(s)}(y)\big) L_{1}( F^{(s)}_{1}(x))
d F^{(s)}(x,y)
\Bigg) .
\e*
\end{theorem}
To normalize the test, we consider the following estimator
\begin{align*}
\w \sigma^2{(1, 2)} & =
\di\frac{(1-a_{1, 2})}{n_{1}}
\di\sum_{i=1}^{n_{1}}(M_{i}^{(1)} - \overline{M}^{(1)})^2
+\di\frac{a_{1, 2}}{n_{2}}
\di\sum_{i=1}^{n_{2}}(M_{i}^{(2)} - \overline{M}^{(2)})^2,
\end{align*}
 \text{ with }
 \begin{align*}
 \overline{M}^{(s)}  = \di\frac{1}{n_s}\sum_{i=1}^{n_s}M_{i}^{(s)},
 \ \
 \text{for \ } s=1,2,
\end{align*}
where
\be*
M_{i}^{(s)} = L_{1}(\w U^{(s)}_{i,1})L_{1}(\w U^{(s)}_{i,2})
&+&
\di\frac{2\sqrt{3}}{n_s} \di\sum_{k=1}^{n_s}\Big( \mathds{1}\big(X^{(s)}_{i,1}\leq X^{(s)}_{k,1}\big)-\w U^{(s)}_{k,1}\Big) L_{1}(\w U^{(s)}_{k,2})
\\
&+&
\di\frac{2\sqrt{3}}{n_s} \di\sum_{k=1}^{n_s}\Big( \mathds{1}\big(X^{(s)}_{i,2}\leq X^{(s)}_{k,2} \big)-\w U^{(s)}_{k,2}\Big) L_{1}(\w U^{(s)}_{k,1}).
\e*
\begin{proposition}\label{prop1_ind}
 Under $ H_0$,
\begin{align*}
\w\sigma^2(1, 2) & \overset{\P}{\longrightarrow}
\sigma^2(1, 2).
\end{align*}
\end{proposition}
We then deduce the limit distribution under the null.
\begin{corollary} \label{cor1_ind}
Assume that  {\bf (A)} holds. Under $ H_0$,
$V_{D(n)}^{(1, 2)}/\w \sigma^2(1, 2)$ converges in law towards  a chi-squared distribution $\chi_1^2$  as $n_{1}, n_2  \rightarrow +\infty$.
\end{corollary}

\section{$K$-sample case}

We restrict our attention to 
the iid case here. 
The paired case 
is treated in Appendix \ref{appendIndep}. 

Our aim is to generalize the two-sample case by considering a series of embedded statistics. Each new statistic will include a new pair of populations to be compared. We will use the first rule  (\ref{rule})  to select a potentially different copula coefficient between each pair.
A  second rule will then be considered to select a possibly different pair between all populations.  
To select the pairs of populations  we introduce the following set of indices:
\be*
{\cal V}(K) & = &
\{(\ell,m)\in \N^2 ; 1\leq \ell<m \leq K\}.
\e*
Clearly, ${\cal V}(K) $ contains   $v(K)=K(K-1)/2$ elements which represent all the pairs  of populations that we want to compare and   that can be ordered as follows:
 we write $(\ell,m) <_{\cal V} (\ell',m')$ if $\ell<\ell'$, or $\ell=\ell'$ and $m<m'$, and we denote by $r_{\cal V}(\ell,m)$ the associated rank of $(\ell,m)$ in ${\cal V}(K)$.
 This can be seen as a natural order (left to right and top to bottom) of the elements of the upper triangle of a $K\times K$ matrix as represented below:
 \be*
 \begin{array}{c c c c c  }
 (1,2)  & (1,3)  &\ldots & \ldots & (1,K)
 \\
 & (2,3)   &\ldots & \ldots & (2,K)
 \\
 & & \ddots & &
 \\
 & & & & (K-1, K)
 \end{array}
 \e*
 We see at once that $r_{\cal V}(1,2) = 1 , r_{\cal V}(1,3)=2$ and more generally, for $\ell, m \in {\cal V}(K) $ we have
 \be*
 r_{\cal V}(\ell,m)=K(\ell-1)-\frac{\ell(\ell+1)}{2}+m.
  \e*
We construct an embedded series of statistics as follows:
\be*
V_1 =
V^{(1,2)}_{D(\m n)}, \ \
V_2  =  V^{(1,2)}_{D(\m n)}+V^{(1,3)}_{D(\m n)}, \ \
\ldots, \  \
V_{v(K)}  = V^{(1,2)}_{D(\m n)} + \cdots + V^{(K-1,K)}_{D(\m n)},
\e*
or equivalently,
\begin{eqnarray*}
V_{k}
& = &
\di\sum_{(\ell,m) \in {\cal V}(K); r_{\cal V}(\ell,m)\leq k} V^{(\ell,m)}_{D(\m n)},
 \e*
 where $D(\m n)$ is given by (\ref{rule}) and $V^{(\ell,m)}_{D(\m n)}$ is defined as in  (\ref{sumV}), replacing the pair index  $(1,2)$ by $(\ell,m)$.
We have $V_1 < V_2 < \ldots < V_{v(K)}$. The first statistic $V_1$ compares the first two populations 1 and 2. The second statistic $V_2$ compares the populations 1 and 2, and, in addition, the populations 1 and 3. And more generally, the statistic $V_k$ compares $k$ pairs of populations. For each $1<k<v(K)$, there exists a unique pair $(\ell,m)$ such that
$r_{\cal V}(\ell,m)=k$.
To choose automatically the appropriate number of pairs $k$  we introduce the following penalization procedure, mimicking the Schwarz criterion procedure \cite{schwarz2}:

\begin{equation} \label{rule3}
s(\m n)= \min \Big\{\argmax_{1 \leq k \leq v(K)} \big({V}_{k} - k  q_{\m n}\big)  \Big\},
\end{equation}
where  $q_{\m n} $ is a penalty term.   
The choice of $q_{\m n}$ is discussed in Remark \ref{remark3}.
We will need the following assumption:
\begin{description}
\item{{\bf (A')}}
$d(n_i)^{(p+5)} = o(q_{\m n}), {\rm \ for \ } i=1,\ldots, K.$
\end{description} 
The following result shows that, under the null, the penalty will choose the first element of ${\cal V}(K)$ asymptotically. This means that all other pairs are not significantly different under the null and do not contribute to the statistic. 

\begin{theorem}
\label{thm3_ind}
Assume that {\bf (A)} and  {\bf (A')}  hold.  Then under $ H_0$,
$s(\m n)$ converges in probability towards 1 as $n_1, \ldots, n_K \rightarrow +\infty$.
\end{theorem}
\begin{corollary}\label{cor2_ind}
Assume that {\bf (A)} and {\bf (A')} hold.
Then under $ H_0$,
$V_{s(\m n)}/\w \sigma^2{(1,2)}$ converges in law towards a $\chi^2_1$ distribution as $n_1, \ldots, n_K \rightarrow +\infty$.
\end{corollary}
Then the final data-driven test statistic is given by
\be*
V & = &
V_{s(\m n)}/\w \sigma^2{(1,2)}.
\e*

\begin{remark}\label{remark3}
In the classical smooth test approach (see \cite{ledwina1994data}), the standard penalty in the univariate case is $q_n=p_n= \log(n)$, a choice closely linked to the Schwarz criteria \cite{schwarz2} as detailed in \cite{KL95}. Here, we extend this approach to the multivariate case with the following generalization:
\begin{equation}
    q_{\bm{n}} = p_{\bm{n}} = \alpha \log\left(\frac{K^{(K-1)}n_1 \cdots n_K}{(n_1+\cdots+n_K)^{K-1}}\right).
\end{equation}
Proposition \ref{thm4} demonstrates that this choice is sufficient for detecting alternatives. In practical applications, the introduction of the factor $\alpha$ serves to stabilize the empirical level, bringing it closer to the asymptotic one. Details on the automatic selection of $\alpha$ can be found in Appendix \ref{annexTuning}, offering a straightforward calibration of the test.

It's worth noting that in \cite{inglot}, a comparison between this Schwarz penalty and the Akaike penalty was conducted. The latter proposes a constant value for $p_{\bm{n}}$ or $q_{\bm{n}}$, providing an alternative approach to calibrating the test.

Finally, in the paired case where $n:=n_1=\ldots = n_K$, we opt for $q_n=p_n= \alpha \log(n)$.
\end{remark}

 \section{Alternative hypotheses}
 \label{secalternatives}
We  consider the following series of  alternative hypotheses: 
for $k \in \{1, \ldots,\\ v(K)\}$
\be*
H_1(k)
\!:\!\left\{
\begin{array}{ll} \!{\rm \! if \ } r_{\cal V}(\ell,m) < k, C_{\ell} {\rm \ and \ } C_{m} {\rm \ have \ the \ same  \ copula \ coefficients} \\
\! {\rm \! if \ } r_{\cal V}(\ell,m) = k, C_{\ell} {\rm \ and \ }C_{m} {\rm \ have \ at \ least\ a \ different  \ copula \ coefficient}.
 \end{array}
 \right.
\e*
The hypothesis $H_1(k)$ asserts that for a given $k$, the populations indexed by $\ell$ and $m$ with $r_{\cal V}(\ell,m)=k$ are the first to exhibit a difference, as per the order defined on ${\cal V}(K)$).
If $k=1$, 
it means
 that the two first copulas $ C_1$   and  $C_2$  have at least one different copula coefficient.
We will need the following assumption:

\begin{description}
\item{{\bf (B)}}
$p_{\m n} = o(\m n)$.
\end{description}
\begin{proposition}\label{prop3_ind}
Assume that {\bf (A)-(A')-(B)} hold. Then under $H_1(k)$, $s(\m n)$ converges in probability towards $k$, as $ n_1, \ldots, n_K \rightarrow +\infty$, and
$V$ 
converges to $+\infty$, that is, $\P(V < \epsilon)\rightarrow 0$, for all $\epsilon >0$.
\end{proposition}
Thus a value of $s(\m n)$  equal to $k$ indicates that the first pairs of populations are equal and that a difference appears from the $k$th pair (following the order on ${\cal V}(K)$).

\section{Numerical study of the test}

We  choose the penalty  $q_{\m n}=p_{\m n} = \alpha \log(K^{(K-1)}n_1 \ldots n_K/(n_1+\ldots+n_K)^{K-1})$, as indicated in Remark \ref{remark3}.  
In our proofs, we set $\alpha=1$  for simplicity. 
However, in practice, we enhance this tuning factor empirically using the data-driven procedure outlined in Appendix \ref{annexTuning}.

Concerning the value of $d(\m n)$, conditions {\bf (A)} and {\bf (A')} are asymptotic conditions and 
from our experience setting $d(\m n)=3$ or $4$ is enough to have a very fast procedure which detects alternatives where 
 copulas differ by a coefficient with a norm less than or equal to $d(\m n)$. This parameter can be modified in the package 'Kcop'. In our simulation, we fixed $d(\m n)=3$. 
 The
nominal level is equal to  $\alpha=5\%$.


\subsection{Simulation design}

We consider the following  copula families:  Gaussian,  Student, Gumbel,  Frank,
Clayton, and Joe Copulas (briefly denoted by Gaus, Stud, Gumb, Fran, Clay and Joe). 
For the explicit  forms and properties
of these copulas, we refer the reader to
\cite{Nelsen}.
For each copula $C$, the sample is generated  with a given
Kendall's $\tau$ parameter, and we denote it briefly by $C(\tau)$.
When $\tau$ is close to zero the variables are close to the independence. Conversely, if $\tau $ is close to 1 the dependence becomes linear.

In our simulation, we compute empirical levels and empirical powers as the percentage of rejections 
under the null and alternative hypotheses based on $1000$ replicates. We consider the following scenarios: 
\bit
\item 
We first consider the two-sample case where we compare our test procedure to that proposed in \cite{remillard} which is the competitor we found for dependent as well as independent bivariate observations. Both methods give very similar results. 
\item 
Then, we consider two cases:
 a $5$-sample
case and  a $10$-sample case. In both situations, alternatives are constructed by modifying  $\tau$. 

\item 
We also compare the performance of the smooth test to the approach developed in  \cite{quessy16} in the $K$-sample case, with $K=2,3,4$,  restricting our study to sub-samples from the observations as done in \cite{quessy16,quessy21}.
\item 
A $6$-population case is studied where we change copulas, keeping the same $\tau$. 
\item Finally an additional simulation study is proposed in Appendix \ref{appendstudent}. We compared three Student copulas with $df=5$ and with $\tau =0.4$ or $0.6$. 
\eit

\subsection{Simulation results in the two-sample case}

In this case ($K=2$) we consider the procedure of
\cite{remillard} as a competitor. Let us recall that
this approach is based on the Cramer-von-Mises
statistic between the two empirical copulas and
an approximate p-value is obtained through
the multiplier technique with $1000$ replications.
They also proposed a R package denoted by {\it Twocop}. 
By extension, we call our  
 R package {\it Kcop}.

Here we fix the dimension $p=2$. The following groups
of scenarios are  considered:
\begin{enumerate}
\item $\mathcal{A}2:\textbf{50-50}$: it includes six alternatives of size $n_1=n_2=50$ which are:
\begin{itemize}
    \item $A2norm$: $C_1=Gaus(\tau_1=0.2)$ and $C_2=Gaus(\tau_2\in \{0.1,0.2,\ldots,\\0.9\})$
     \item $A2stu$: $C_1=Stud(df=17,\tau_1=0.2)$ and $C_2=Stud(df=17,\tau_2\in \{0.1,0.2,\ldots,0.9\})$ where $df$ is a degree of freedom
      \item $A2gum$: $C_1=Gumb(\tau_1=0.2)$ and $C_2=Gumb(\tau_2\in \{0.1,0.2,\ldots,\\0.9\})$
       \item $A2fran$: $C_1=Fran(\tau_1=0.2)$ and $C_2=Fran(\tau_2\in \{0.1,0.2,\ldots,\\0.9\})$
        \item $A2clay$: $C_1=Clay(\tau_1=0.2)$ and $C_2=Clay(\tau_2\in \{0.1,0.2,\ldots,0.9\})$
         \item $A2joe$: $C_1=Joe(\tau_1=0.2)$ and $C_2=Joe(\tau_2\in \{0.1,0.2,\ldots,0.9\})$
\end{itemize}
\item $\mathcal{A}2:\textbf{50-100}=\mathcal{A}2:\textbf{50-50}$ with $n_1=50$ and $n_2=100$
\item $\mathcal{A}2:\textbf{100-50}=\mathcal{A}2:\textbf{50-50}$ with $n_1=100$ and $n_2=50$
\item $\mathcal{A}2:\textbf{100-100}=\mathcal{A}2:\textbf{50-50}$ with $n_1=100$ and $n_2=100$
\end{enumerate}
Recall that this methodology to evaluate  the finite sample
performance was 
proposed in \cite{remillard}. 
We follow their designs with the same sample sizes $(n_1,n_2)\in
\big\{(50,50),(50,100),(100,50),(100,100) \big\}$. 
Such scenarios coincide with the null hypothesis when $\tau_2=0.2$.

The results are very similar for all scenarios and we present  the $A2norm$ alternatives in this section,  reserving the remaining  results for  Appendix \ref{appendRemillard}. 
Figures \ref{double1}-\ref{double2} illustrate  that both methods ({\it Twocop}  and {\it Kcop})  exhibit highly comparable performance.  As expected, the more different the Kendall tau, the greater the power.  In our simulation, the tau associated with $C_1$ is fixed and equal to 0.2. The tau associated with $C_2$ varies and the power is maximal (100\%) when it is greater than or equal to  0.7. Conversely, the power is  minimal (approaching 5\%) when the tau is set at 0.2, corresponding to the null hypothesis.

\FloatBarrier

\begin{figure}[!htbp]
    \centering 
\begin{subfigure}{0.49\textwidth}
  \includegraphics[width=\linewidth]{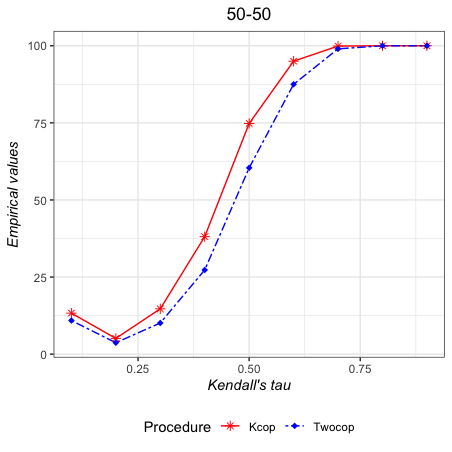}
\end{subfigure}\hfil 
\begin{subfigure}{0.49\textwidth}
  \includegraphics[width=\linewidth]{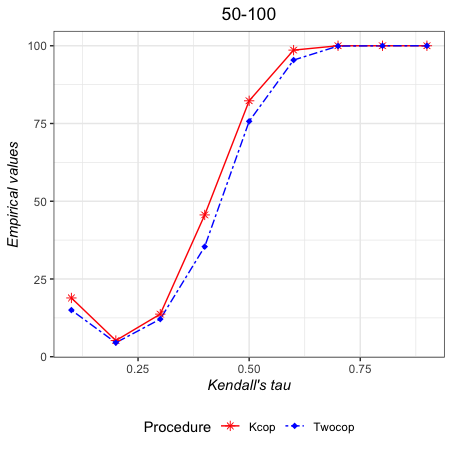}
\end{subfigure}\hfil 
\caption{Two-sample case: \% of rejections under $\mathcal{A}2norm:\textbf{50-50}$ (left) and $\textbf{50-100}$ (right).}
\label{double1}
\end{figure}

\FloatBarrier
\begin{figure}[!htbp]
    \centering 
\begin{subfigure}{0.49\textwidth}
  \includegraphics[width=\linewidth]{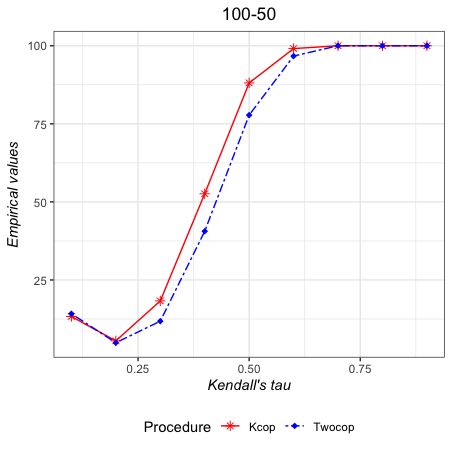}
\end{subfigure}\hfil 
\begin{subfigure}{0.49\textwidth}
  \includegraphics[width=\linewidth]{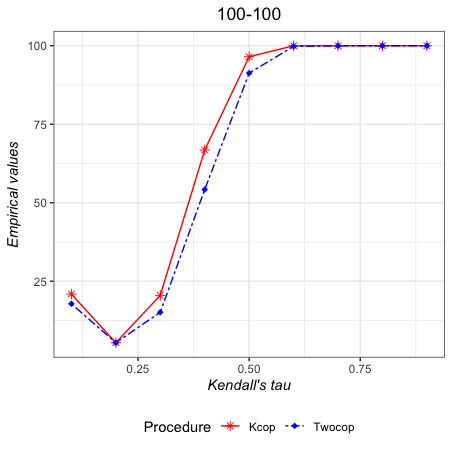}
\end{subfigure}\hfil 
\caption{Two-sample case: \% of rejections  under $\mathcal{A}2norm:\textbf{100-50}$ (left) and $\textbf{100-100}$ (right).}
\label{double2}
\end{figure}

\subsection{Five-sample case}

In this case ($K=5$) we fix  $p=3$
and we consider the same size for all samples, that is
$n=n_1=n_2=n_3=n_4=n_5$ $\in
\{50,100,200,\ldots,900,1000\}$.
We fixed a theoretical level
$\alpha=5\%$.

~\\{\bf Null hypotheses}: under the null  hypothesis we consider the same copulas (Gaussian, Student with degree of freedom = 17, Gumbel, Frank, Clayton, Joe)  with
    three levels of dependence: $\tau=0.1$ (low dependence), $\tau=0.5$ (middle dependence) and $\tau=0.8$ (high dependence).

~\\
{\bf Alternatives with different tau}: we consider  the following  alternatives hypotheses with  $C_1, \ldots, C_5$ in the same copula family but with different $\tau$ as follows  
\begin{itemize}
    \item \textbf{Alt1}:  $C_1(0.3)=C_2(0.3)=C_3(0.3)=C_4(0.3)$ and $C_5(0.1)$ 
    \item \textbf{Alt2}: $ C_1(0.1)$ and $C_2(0.55)=C_3(0.55)= C_5(0.55)$, and $C_4(0.3)$
        \item \textbf{Alt3}: $C_1(0.1)$ and $C_2(0.8)=C_3(0.8)=C_5(0.8)$, and $C_4(0.3)$
\end{itemize}

{\bf Alt1} contains only one different population. Concerning {\bf Alt2} and {\bf Alt3}, they differ solely in their Kendall's tau, allowing us to highlight its effect.

Table \ref{level5} presents empirical levels (in $\%$) with respect to sample sizes when $\tau=0.1, 0.5$ and $0.8$, respectively. In each case, one can observe that the empirical level is close to the theoretical 5\% as soon as $n$ is greater than 200. For $n=50$ or $100$, two phenomena emerge: the empirical level appears larger than the theoretical level when $\tau $ is small  
and smaller than the theoretical level when $\tau$  is large.  
Hence, with fewer observations, the procedure more readily identifies identical copulas when their dependence structure is stronger. This leads to the following recommendations: for a small size ($n<200$) if the estimation of $\tau$ is close to 0.1, it is advisable to adopt a more conservative approach (choosing a larger theoretical level, e.g., around 0.09). 
Conversely, if the estimation of $\tau$ is close to 0.9, it is preferable to be anticonservative (choosing a lower theoretical level around  0.02). This implies a slight reduction in power in the first case, while power increases in the second case. 
A tuning procedure could be considered, incorporating a data-driven criterion based on the estimation of $\tau$.

\begin{table*}[ht]
\caption{Empirical levels (in \%) for the five-sample test.}
\label{level5}
\centering
\begin{tabular}{ccccccc}
\hline
& \multicolumn{6}{c}{Models} \\ \cline{2-7}
$n$  & Gaussian & Student  & Gumbel & Frank  & Clayton & Joe \\
  \hline 
& \multicolumn{6}{c}{Kendall tau $\tau=0.1$} \\ \cline{2-7}
50 & 11.4 & 10.5 & 10.0 & 11.1 & 10.3 & 11.4 \\
100 & 10.0 & 8.4 & 8.1 & 7.6 & 8.1 & 9.1 \\
200 & 7.6 & 8.0 & 6.2 & 6.3 & 5.8 & 7.4 \\
300 & 6.9 & 7.3 & 6.6 & 7.5 & 6.5 & 6.3 \\
400 & 6.4 & 5.7 & 7.1 & 4.7 & 5.7 & 7.4 \\
500 & 5.1 & 4.8 & 4.9 & 7.0 & 5.9 & 5.5 \\
600 & 5.6 & 6.5 & 4.4 & 5.1 & 5.1 & 6.0 \\
700 & 5.0 & 5.3 & 6.1 & 5.5 & 4.4 & 6.5 \\
800 & 5.1 & 6.8 & 4.8 & 5.5 & 5.4 & 6.0 \\
900 & 5.6 & 6.2 & 6.3 & 5.7 & 6.5 & 6.8 \\
1000 & 5.9 & 5.5 & 6.0 & 5.3 & 5.2 & 5.0 \\
\cline{2-7}
& \multicolumn{6}{c}{Kendall tau $\tau=0.5$} \\ \cline{2-7}
50 & 5.4 & 4.0 & 4.0 & 4.1 & 5.0 & 3.2 \\
100 & 6.0 & 3.7 & 5.0 & 5.4 & 5.1 & 2.8 \\
200 & 4.9 & 5.0 & 5.6 & 5.5 & 6.0 & 4.8 \\
300 & 5.7 & 3.9 & 4.6 & 4.9 & 5.6 & 4.0 \\
400 & 4.7 & 3.9 & 4.9 & 5.1 & 4.6 & 4.6 \\
500 & 4.4 & 3.6 & 3.6 & 5.5 & 4.4 & 4.5 \\
600 & 4.8 & 5.0 & 3.2 & 4.2 & 4.7 & 5.5 \\
700 & 5.4 & 5.5 & 5.0 & 6.0 & 5.0 & 4.6 \\
800 & 4.9 & 4.6 & 4.6 & 3.7 & 4.5 & 4.4 \\
900 & 4.6 & 5.0 & 4.2 & 6.1 & 4.2 & 4.0 \\
1000 & 4.2 & 4.6 & 4.1 & 4.9 & 5.8 & 3.5 \\

\cline{2-7}
& \multicolumn{6}{c}{Kendall tau $\tau=0.8$} \\ \cline{2-7}
      50 & 1.0 & 0.6 & 0.6 & 0.7 & 3.0 & 0.4 \\
     100 & 2.6 & 1.9 & 2.2 & 2.9 & 4.5 & 1.4 \\
      200 & 4.1 & 3.1 & 3.5 & 3.9 & 5.3 & 3.0 \\
     300 & 4.0 & 3.3 & 4.5 & 3.4 & 5.4 & 2.1 \\
    400 & 3.5 & 3.4 & 4.3 & 4.2 & 5.5 & 3.9 \\
      500 & 4.9 & 3.9 & 3.4 & 3.8 & 4.0 & 3.6 \\
    600 & 4.6 & 3.9 & 4.1 & 4.5 & 5.1 & 4.8 \\
      700 & 4.0 & 5.4 & 4.0 & 4.4 & 5.8 & 3.7 \\
    800 & 4.5 & 4.6 & 5.0 & 5.0 & 4.8 & 4.1 \\
    900 & 4.4 & 4.1 & 3.8 & 5.1 & 4.8 & 4.2 \\
   1000 & 3.7 & 5.4 & 3.8 & 5.6 & 5.4 & 4.1 \\

\hline
\end{tabular}
\end{table*}

   Concerning the empirical power, Tables \ref{tab1}-\ref{tab3} contain all results under the alternatives. We 
    omit some large sample size results where empirical powers are equal to 100\%.
   It is important to note that, even for a sample size equal to $1000$, the program runs very fast.
   It can be seen for alternatives {\bf Alt2} and {\bf Alt3} that the empirical powers are extremely high even for small sample sizes.
   The first series of alternatives yields lower empirical powers since only one copula differs with a slight change in $\tau$.

\FloatBarrier


\begin{table*}[ht]
\caption{Empirical powers (in \%)  under alternative  {\bf Alt1} (five-sample case).
}\label{tab1}
\centering
\begin{tabular}{ccccccc}
\hline
& \multicolumn{6}{c}{Alternatives} \\ \cline{2-7}
  & Gaussian & Student  & Gumbel & Frank  & Clayton & Joe \\
{{$n=50$ }}
  & 39.9 & 35.7 & 35.6 & 36.6 & 35.9 & 35.5  \\
{{$n=100$}}
 & 64.1 & 61.8 & 60.3 & 64.0 & 61.1 & 60.7 \\
{{$n=200$}}
  & 91.5 & 88.4 & 87.5 & 91.1 & 89.9 & 87.7  \\
{{$n=300$ }}
  & 97.9 & 98.0 & 97.7 & 98.2 & 97.3 & 97.2 \\
{{$n=400$ }}
  & 99.8 & 99.7 & 99.6 & 99.8 & 99.7 & 99.8  \\
{{$n=500$ }}
 & 100 & 100 & 100 & 100 & 100 & 99.9  \\
{{$n=600$}}
  & 100 & 100 & 100 & 100 & 100 & 100  \\
\hline
\end{tabular}
\end{table*}

\FloatBarrier

\begin{table*}[ht]
\caption{\label{tab2}  Empirical powers (in \%)   under alternative  {\bf Alt2} (five-sample case).
}
\centering
\begin{tabular}{ccccccc}
\hline
  & \multicolumn{6}{c}{Alternatives} \\ \cline{2-7}
& Gaussian & Student  & Gumbel & Frank  & Clayton & Joe  \\
{{$n=50$ }}
  & 97.8 & 97.6 & 96.3 & 98.6 & 97.4 & 95.6  \\
{{$n=100$}}
  & 100 & 100 & 99.9 & 100 & 100 & 100 \\
{{$n=200$}}
& 100 & 100 & 100 & 100 & 100 & 100  \\
\hline
\end{tabular}
\end{table*}

\FloatBarrier

\begin{table*}[ht]
\caption{\label{tab3} Empirical powers (in \%)  under alternative  {\bf Alt3} (five-sample case).
}
\centering
\begin{tabular}{ccccccc}
\hline
  & \multicolumn{6}{c}{Alternatives} \\ \cline{2-7}
  & Gaussian & Student  & Gumbel & Frank  & Clayton & Joe   \\
{{$n=50$ }}
 & 100 & 100 & 100 & 100 & 100 & 100  \\
\hline
\end{tabular}
\end{table*}

\subsection{Ten-sample case}

Analogously to the previous $5$-sample case,
we consider
 null hypotheses with Gaussian, Student, Gumbel, Frank, Clayton, and Joe copulas. We fixed $p=2$. 
We consider the following alternatives where only one copula differs from the others.
\bit
\item
$\textbf{Alt4}$:
$C_1(0.1)= C_2(0.1) =
 \cdots=C_{9}(0.1)$ and $C_{10}(0.55)$
\eit
Empirical levels seem to tend fast to 0.5 and are relegated in  Appendix \ref{appendempirical}.
Table \ref{tab4} shows empirical powers under alternatives {\bf Alt4}. We only treat the cases where $n=50$ and $100$, as beyond these values,  all empirical powers are equal to $100\%$.  Remarkably, even for such small sample sizes, we observe very good behavior of the test even with small sample sizes.


\FloatBarrier

\begin{table*}[ht]
\caption{\label{tab4} Percentage of rejection under alternative {\bf Alt4}  (ten-sample case).
}
\centering
\begin{tabular}{ccccccc}
\hline
  & \multicolumn{6}{c}{Alternatives} \\ \cline{2-7}
  & Gaussian & Student  & Gumbel & Frank  & Clayton & Joe   \\
{{$n=50$ }}
 & 98.0 & 96.7 & 96.2 & 97.9 & 97.1 & 97.3  \\
{{$n=100$}}
  & 100 & 100 & 100 & 100 & 100 & 100 \\
\hline
\end{tabular}
\end{table*}

\FloatBarrier

\subsection{Alternatives with the same Kendall's tau}
We consider a last  alternative hypothesis with  $C_1, \ldots, C_6$ %
which are the six  copulas defined in the null hypothesis models above
%
all with the same $\tau=0.55$ and with a dimension $p$ up to 5 as follows  
\begin{itemize}
    \item \textbf{Alt5}:  $\tau=0.55$;  $C_1=Gauss$, $C_2=Student$, $C_3=Gumbel$, $C_4=Frank$, $C_5=Clayton$, $C_6=Joe$, and the dimension $p\in\{2,3,4,5\}$.

\end{itemize}

Empirical powers are presented in Table \ref{tab5}. It can be seen that the power increases with the dimension $p$ when the sample size is less than $n=300$: it is then easier to detect differences between the dependence structures of the vectors.
When $n\geq 300$, the empirical power is stable and equal to $100\%$ in all scenarios.

\FloatBarrier

\begin{table*}[ht]
\caption{\label{tab5} Empirical powers (in \%) under alternative {\bf Alt5} (6-sample case).
}
\centering
\begin{tabular}{ccccc}
\hline
 Dimension  & $p=2$ & $p=3$  & $p=4$ & $p=5$   \\
 \hline 
{{$n=50$ }}   & 1.2 &  3.1 & 14.8  & 20.0 \\\hline
 {{$n=100$ }} & 2.0  & 27.3  & 73.6    & 79.1 \\\hline
 {{$n=200$ }} & 19.8 & 89.9  & 99.8   & 100   \\\hline
{{$n=300$ }} & 60.3 & 100.0 & 100.0 & 100    \\\hline
 {{$n=400$ }} & 90.9 & 100.0 & 100.0 & 100    \\\hline
 {{$n=500$ }} & 98.3 & 100.0 & 100.0 & 100    \\\hline
\end{tabular}
\end{table*}

\FloatBarrier

\subsection{Testing the equality of all the bivariate sub-copulas of copulas}
The purpose of this section is to compare the performance of our test with that obtained by  \cite{quessy16}.
We follow the same design  (see Tables 3 in \cite{quessy16}) and we adopt the same notation. More precisely, we simulated data $\bm U=(U_1,\ldots,U_{2K})\sim C$, where $C$ is a $2K$-dimensional copula and we examine the equality of all the bivariate sub-copulas of $\bm U$, that is
\be*
\di C^{(U_1,U_2)}=C^{(U_3,U_4)}=\cdots=C^{(U_{2K-1},U_{2K})}.
\e*
We denote by $N(\theta$) the model where  $\bm U$ is generated by the $2K$-variate normal copula and by $T(\theta)$ the model where
 $\bm U$ is generated by the Student copula with $\nu=3$ degrees of freedom, where the correlation matrix
 $\Sigma$ is such that $\theta=\Sigma_{1,2}=\Sigma_{2,1}$ and $\Sigma_{i,j}=0.2$ for all $(i,j)\neq \{(1,2);(2,1)\} $

 We compare our procedure ($Kcop$) to the following quadratic functional procedures proposed in  \cite{quessy16}:
 \begin{itemize}
     \item Cramér-von Mises ($CvM$) statistic,
     \item 
     Two characteristic function statistics, denoted as ($Cf_1$,$Cf_2$), correspond to the weights functions of normal and double-exponential distributions, respectively
     \item Diagonal statistics (Dia).
 \end{itemize}
We refer the reader to \cite{quessy16} for more detail and to \href{www.uqtr.ca/MyMatlabWebpage}{code} for the program.

The results are provided in Tables \ref{TabQuessy0} and \ref{TabQuessy}. There is no overarching conclusion that allows determining a superior method. The various statistics seem to yield fairly similar results, except in the case of $K=4$, where the emprical powers  associated with our  test statistic appear to be generally superior.

\FloatBarrier

\begin{table}[!h]
\caption{ Empirical levels  for   different models studied in \cite{quessy16} with  sample size $n=50$ and $n=100$.
}\label{TabQuessy0}
\centering
\begin{adjustbox}{max width=1\textwidth,center}  
\begin{tabular}{ccccccccccccc}
\hline
& & \multicolumn{5}{c}{$n=50$} && \multicolumn{5}{c}{$n=100$}
\\
& & \hhline{----- ~ -----}
 {$K$} &
{Approaches} & \multicolumn{1}{c}{$CvM$} & \multicolumn{1}{c}{$Dia$} & \multicolumn{1}{c}{$Cf_1$} & \multicolumn{1}{c}{$Cf_2$} & \multicolumn{1}{c}{$Kcop$} & &  \multicolumn{1}{c}{$CvM$} & \multicolumn{1}{c}{$Dia$} & \multicolumn{1}{c}{$Cf_1$} & \multicolumn{1}{c}{$Cf_2$} & \multicolumn{1}{c}{$Kcop$}
\\
& & \hhline{----- ~ -----}
%
%
\multirow{2}{*}{$2$} & $N(.2)$ & 4.7 & 5.2 & 6.2 & 6.1 & 6.0 && 4.7 & 3.9 & 4.7 & 4.7 & 4.0 \\  
 & $T(.2)$ & 4.1 & 3.9 & 4.8 & 4.6 & 6.0 & &4.4 & 4.6 & 5.2 & 5.3 & 5.0
 \\
& & \hhline{----- ~ -----}
%
%
\multirow{2}{*}{$3$} & $N(.2)$ & 3.3 & 4.8 & 5.7 & 4.8 & 4.0 && 2.9 & 4.1 & 3.3 & 3.8 & 4.0 \\  
 & $T(.2)$ & 3.0 & 4.4 & 3.9 & 3.7 & 6.0 && 4.0 & 5.1 & 4.3 & 4.5 & 6.0 \\
& & \hhline{----- ~ -----}
%
%
\multirow{2}{*}{$4$} & $N(.2)$ & 3.4 & 4.1 & 5.7 & 4.9 & 5.0 && 2.3 & 3.0 & 3.5 & 3.3 & 6.0 \\  
 & $T(.2)$ & 1.7 & 4.9 & 3.5 & 3.0 & 6.0 && 4.4 & 4.5 & 6.2 & 6.1 & 6.0 \\
\hline
%
\end{tabular}
\end{adjustbox}
\end{table}
\begin{table}
\caption{\label{TabQuessy} Empirical powers for  different alternatives studied in \cite{quessy16} with   $n=50$ and $n=100$.
}
    \centering
\begin{adjustbox}{max width=1\textwidth,center}  
\begin{tabular}{ccccccccccccc}
\hline
%
&&  \multicolumn{5}{c}{$n=50$} & &\multicolumn{5}{c}{$n=100$}  \\
& & \hhline{----- ~ -----}
{$K$} & 
{Approaches} & \multicolumn{1}{c}{$CvM$} & \multicolumn{1}{c}{$Dia$} & \multicolumn{1}{c}{$Cf_1$} & \multicolumn{1}{c}{$Cf_2$} & \multicolumn{1}{c}{$Kcop$} & &  \multicolumn{1}{c}{$CvM$} & \multicolumn{1}{c}{$Dia$} & \multicolumn{1}{c}{$Cf_1$} & \multicolumn{1}{c}{$Cf_2$} & \multicolumn{1}{c}{$Kcop$} \\
%
%
& & \hhline{----- ~ -----}
\multirow{4}{*}{$2$} & $N(.4)$ & 16.8 & 15.5 & 21.3 & 19.8 & $\bm{32.0}$ & & 22.5 & 21.8 & 28.5 & 26.4 & $\bm{43.0}$ \\  
 & $T(.4)$ & 44.2 & 48.1 & $\bm{48.9}$ & 48.1 & 35.0 && 78.0 & 75.4 & $\bm{82.2}$ & 81.2 & 40.0 \\
 & $N(.6)$ & 51.3 & 48.9 & 62.2 & 58.1 & $\bm{68.0}$ & & 84.1 & 81.7 & 90.1 & 87.7 & $\bm{94.0}$ \\
 & $T(.6)$ & 88.6 & 92.2 & $\bm{90.8}$ & 90.7 & 60.0 & & $\bm{99.9}$ & $\bm{99.9}$ & $\bm{99.9}$ & $\bm{99.9}$ & 84.0 \\
& & \hhline{----- ~ -----}
%
%
\multirow{4}{*}{$3$} & $N(.4)$ & 12.0 & 12.7 & 17.8 & 16.1 & $\bm{64.0}$ & & 20.3 & 21.1 & 25.2 & 23.6 & $\bm{74.0}$ \\  
 & $T(.4)$ & 40.7 & 47.7 & 47.3 & 45.7 & $\bm{70.0}$ && 76.7 & 76.6 & $\bm{81.2}$ & 79.4 & 73.0 \\
 & $N(.6)$ & 47.3 & 48.2 & 63.0 & 56.8 & $\bm{92.0}$ && 85.8 & 84.7 & 91.5 & 88.6 & $\bm{99.0}$ \\
 & $T(.6)$ & 90.1 & $\bm{93.5}$ & 93.2 & 92.4 & 89.0 && 100.0 & 99.8 & 99.9 & $\bm{100.0}$ & 97.0 \\
& & \hhline{----- ~ -----}
\multirow{4}{*}{$4$} & $N(.4)$ & 9.8 & 12.0 & 16.1 & 14.1 & $\bm{78.0}$ & &19.6 & 20.6 & 27.1 & 25.0 & $\bm{84.0}$ \\  
 & $T(.4)$ & 34.6 & 41.7 & 43.9 & 42.3 & $\bm{84.0}$ && 74.4 & 75.5 & 78.8 & 77.6 & $\bm{86.0}$ \\
 & $N(.6)$ & 43.4 & 45.8 & 57.2 & 53.5 & $\bm{95.0}$ & & 81.5 & 80.3 & 88.9 & 86.4 & $\bm{100.0}$ \\
 & $T(.6)$ & 86.3 & 91.7 & 91.0 & 90.5 & $\bm{96.0}$ && 99.8 & 99.8 & $\bm{100.0}$ & $\bm{100.0}$ & 99.0 \\
\hline 
\end{tabular}
\end{adjustbox}
\end{table}

\FloatBarrier

\section{Real datasets applications}

\subsection{Biology data}

We analyze Fisher’s well-known  Iris dataset.
The data consists of fifty observations of four measures: Sepal Length ($SL$), Sepal Width ($SW$), Petal Length ($PL$),
and Petal Width ($PW$), for each of three Species: Setosa, Virginica, and Versicolor.
We then have $K=3$  populations, and the dimension is $p=4$.
The lengths and widths for the three species are represented in Appendix \ref{annexeIris}.
In \cite{dhar} the authors
show that multivariate normal distributions seem to fit the data well for all three Iris species. Looking at their mean parameters the 4-dimensional joint distributions seem different but that does not tell us about their dependence structures.

We propose to test the equality of the dependence structure between the four variables $(SL, SW, PL, PW)$ in the three-sample case, that is:
$$
H_0: C_{Setosa}=C_{Virginica}=C_{Versicolor}.
$$
We consider the data as possibly dependent, with the same sample size $n=50$. We then apply the test for paired data.  
We obtain a p-value close to zero ($10^{-11}$)  and a very large test statistic $V=45.9$.
We reject the equality of the dependence structure here. The selected rank  $s( n)$ is equal to 2. It means that the most significant difference is obtained when considering the statistics associated with population 1 versus 2 (Setosa and  Virginica) and population 1 versus 3 (Setosa and Versicolor).

In case of rejection, we can proceed to an ``ANOVA'' type procedure, applying 
 a series of two-sample tests. Table \ref{tabiris}  contains the associated p-values and we conclude with the equality of the dependence structure between  Versicolor and Virginica.

\FloatBarrier

~\\
\begin{table}[!ht]
\caption{P-values for the two-sample test (Iris dataset).}
    \label{tabiris}
    \center
\begin{tabular}{cccc}
\hline
& Setosa & Virginica  & Versicolor
\\
Setosa & 1 & 0.0021 &  $10^{-8}$
\\
Virginica &  0.0021 & 1 & 0.68
\\
Versicolor & $10^{-8}$  & 0.68 & 1
\\
\hline
\end{tabular}
\end{table}

\FloatBarrier

%

\subsection{Insurance data}

Insurance is an area in which understanding  the dependence structure among multiple portfolios is crucial for pricing, especially for risk pooling or price segmentation. To illustrate, we examine the Society of Actuaries Group Medical Insurance Large Claims Database, which contains claims information for each claimant from seven insurers over the period $1997$ to $1999$. 
Each  row in the database presents a summary of
claims for an individual claimant in $27$ fields (columns). 
The first five columns provide general
information about the claimant, the next twelve quantify various types of medical charges and expenses, and the last ten columns summarize details related to the diagnosis. 
For a detailed and thorough description of the data available online, refer to  \cite{soadata},  accessible   on the web page of the 
 \href{https://www.soa.org/resources/experience-studies/2000-2004/research-medical-large-claims-experience-study}{Society of Actuaries}. 
In this context, we focus on  $p=3$
 dimensional variables $\textbf{X} =(X_1,X_2,X_3)$,
 where 
 $X_{1}$ = paid hospital charges, 
$X_{2}$ = paid
 physician charges, 
 $X_{3}$ =  paid other
 charges, 
for all claimants insured by a
 Preferred Provider Organization  plan
 providing exposure for members. This consideration becomes pertinent for risk pooling if the objective is to group together similar charge scenarios or for price segmentation to provide similar guarantees for the charges. We employ a procedure with three scenarios to study the dependence structure of  $\m X$ as follows:

~\\{\bf Three-sample test,   paired case}.
In this case, we consider the same claimants (paired situation) present
over the three periods $1997-1999$. At the end of
the data processing, we obtained three samples of
size $n=6874$ observations.
We analyse the dependence structure of the charges
\textbf{X} between the three years, that is, we test  $H_0: C_{\textbf{X}}^{1997}=C_{\textbf{X}}^{1998}=C_{\textbf{X}}^{1999}$. 
The test concluded with the non-rejection of the equality of the three dependence structures, as evidenced by a p-value = 0.788, 
 a test statistic of $V=0.072$ and a  selected rank equal to $s(n)=1$. Hence, the dependence structure of paid for insured over the three years seems to be similar. It can be an argument for keeping the same distribution of risks on the different charges $X_1, X_2$ and $X_3$.

~\\{\bf Three-sample test,   independent case}. 
Here, we narrow our focus to female claimants. The three populations consist of individuals classified by their relationship with the subscriber, which can be ``Employee'' ($n_E=18144$ observations), ``Spouse'' ($n_S=10969$ observations), or ``Dependent'' ($n_D=10969$ observations), all for the year $1999$. 

Our objective is to test the equality of the dependence structure among the charges $\mathbf{X}$. In this context, we assume independence among the $K=3$ populations. Through our testing procedure, we obtain a p-value close to zero. Consequently, we reject the null hypothesis of equal dependence structure for these charges. 

Subsequently, applying an ANOVA procedure reveals that the two-by-two equalities are rejected for ``Dependent'' vs ``Employee'' and ``Employee'' vs ``Spouse'', with a p-value close to zero in each case. The p-value for ``Dependent'' vs ``Spouse'' is close to one. 

Therefore, the status of being a ``Dependent'' or ``Spouse'' implies a similar dependence structure for the charges, distinct from the status of being an ``Employee''. In the context of risk pooling, differentiating charges between these two groups becomes relevant.

~\\{\bf Ten-sample test,  independent case}.
Here, we analyze data from the year $1999$ where the relationship to the subscriber is ``Employee''. We categorize the charges $\mathbf{X}$ based on age ranges of three years, creating $10$ groups as follows: $G_1=[1936, \: 1938]$, $\ldots$, $G_{10}=[1963, \: 1965]$.

The null hypothesis is $H_0$: the dependence structures of these 10-sample groups are identical. Applying our test procedure, we obtain a p-value close to $0$ and a test statistic of $V=16.20$. Thus, we reject the null hypothesis of equal dependence structure by age at a significant level of $\alpha=5\%$.

There is evidence to suggest that the dependence structure of $\mathbf{X}$ changes over age. We further apply an ANOVA procedure, and the results are presented in Appendix \ref{appendANOVA}, Table \ref{anova_pvalue}, where a two-by-two comparison is proposed. Notably, there are no significant differences between two successive years. Additionally, Group 6 exhibits a similar dependence structure to the other groups, except for Group 3. The disparity increases with the gap between the years, especially between the first age categories and the last ones.

Observing the age range, we identify two clusters: $\{$Group 1, $\ldots$, Group 5$\}$ and $\{$Group 6, $\ldots$, Group 10$\}$. In terms of price segmentation, this allows the formation of two groups with similar dependencies.

\section{Other similar tests}\label{appendOther}
Some extensions of the K-sample test to various null hypotheses have been studied in \cite{quessy16, quessy2017, quessy21}. 
Following his approach we indicate how to adapt the previous test procedure to answer the following hypotheses:
\begin{align*}
    H_0^{RS}: & C^{(\m U^{(1)}, \ldots, \m U^{(K)})} = C^{(1-\m U^{(1)}, \ldots, 1- \m U^{(K)})}
    \\
    H_0^{Exc}: &  C^{(\m U^{(\ell)}, \m U^{(m)})} = C^{\m U^{(m)},  \m U^{(\ell)})} {\rm \ }, \forall \ell\neq m
    \\
   H_0^{ES}: & C^{(\m U^{(1)}, \ldots , \m U^{(K)})} = C^{\m U^{(j_1)}, \ldots  \m U^{(j_K)})}, {\rm \ for \ all \ permutations \ } \m j  {\rm \ of \ }\{1,\ldots, K\}
\end{align*}
Clearly, $H_0^{RS}$ coincides with the radial symmetry, that is $(\m U^{(1)}, \ldots, \m U^{(K)})$ and  $(1-\m U^{(1)}, \ldots, 1- \m U^{(K)})$ have the same joint distribution, while $H_0^{Exc}$ means that copulas are pairwise exchangeable. The exchangeable symmetry is represented by $H_0^{ES}$.
These three hypotheses have been elegantly grouped together and tested in
\cite{quessy16,quessy21}. We can also adapt our procedure to such hypotheses very naturally by considering the density representation given by
(3).
For instance, in the two-sample case, testing $H_0^{RS}$ remains to compare the coefficients $\E\left(L_{\m j_1}(\m U^{(1)})L_{\m j_2}(\m U^{(2)})\right)$ to the coefficients  $\E\left(L_{\m j_1}(1-\m U^{(1)})L_{\m j_2}(1-\m U^{(2)})\right)$ for all $\m j_1, \m j_2$ in $\N^p$.
Asymptotically, under $H_0^{RS}$ the test statistic coincides with the comparison of
$\E\left(L_{1}(U_1^{(1)})L_{1}(U_1^{(2)})\right)$ to   $\E\left(L_{1}(1- U_1^{(1)})L_{1}(1- U_1^{(2)})\right)$ and the selected test statistic is
\begin{align*}
    \di\frac{1}{n}
    \di\sum_{i=1}^{n}
\big(    L_1(\w U^{(1)}_{i,1}) L_1(\w U^{(2)}_{i,1})
    -
    L_1(1-\w U^{(1)}_{i,2}) L_1(1-\w U^{(2)}_{i,2})
\big),
\end{align*}
which has an asymptotic centred normal distribution under  $H_0^{RS}$ with variance similar to that studied in Proposition 1 of the paper.

 In the same way,   $H_0^{Exc}$ consists in comparing 
 $\E\left(L_{\m j_1 }(\m U^{(\ell)})L_{\m j_2}(\m U^{(m)})\right)$ to 
 $\E\left(L_{\m j_1}(\m U^{(m )})L_{\m j_2}( \m U^{(\ell)})\right)$ for all
 $\ell\neq m$.
 Under the null hypothesis, the test statistic coincides simply with the comparison of the first coefficients (the least penalized)
$\E\left(L_{1}(U_1^{(\ell)})L_{1}(U_2^{(\ell)})\right)$ and    
$\E\left(L_{1}(U_1^{(m)})L_{1}(U_2^{(m)})\right)$,  asymptotically.
Then the selected statistic under the null is 
\begin{align*}
    \di\frac{1}{n}
    \di\sum_{i=1}^{n}
\big(    L_1(\w U^{(\ell)}_{i,1}) L_1(\w U^{(\ell)}_{i,2})
    -
    L_1(\w U^{(m)}_{i,1}) L_1(\w U^{(m)}_{i,2})
\big),
\end{align*}
which has asymptotically a centered normal null distribution.

Finally, the same reasoning applies to $H_0^{ES}$ where the test statistic is asymptotically the same as the previous one.

We propose now to compare the performance of our test to the one developed in  \cite{quessy16} for testing the equality of all the bivariate sub-copulas of copulas.
We follow the design given in  \cite{quessy16}(see Table 3). We adopt the same notation and
the same design. More precisely, we simulated data $\bm U=(U_1,\ldots,U_{2K})\sim C$, where $C$ is a $2K$-dimensional copula and we examine the equality of all the bivariate sub-copulas of $\bm U$, that is
\be*
\di C^{(U_1,U_2)}=C^{(U_3,U_4)}=\cdots=C^{(U_{2K-1},U_{2K})}.
\e*
We denote by $N(\theta$) the model where  $\bm U$ is generated by the $2K$-variate normal copula and by $T(\theta)$ the model where
 $\bm U$ is generated by the Student copula with $\nu=3$ degrees of freedom where the correlation matrix
 $\Sigma$ is such that $\theta=\Sigma_{1,2}=\Sigma_{2,1}$ and $\Sigma_{i,j}=0.2$ for all $(i,j)\neq \{(1,2);(2,1)\} $

 We compare our procedure ($Kcop$) to the following quadratic functional 
 procedures proposed in  \cite{quessy16}:
 \begin{itemize}
     \item Cramér-von Mises ($CvM$) statistic,
     \item 
     Two characteristic function statistics, denoted as ($Cf_1$,$Cf_2$), correspond to the weights functions of normal and double-exponential distributions, respectively,
     \item Diagonal statistics (Dia)
 \end{itemize}
We refer the reader to \cite{quessy16} for more details on these procedures and to \href{www.uqtr.ca/MyMatlabWebpage}{code} on their program.

\FloatBarrier

\begin{table}[!h]
\caption{ Empirical levels for different models studied in \cite{quessy16} with  sample size $n=50$ and $n=100$.
}\label{TabQuessy0}
\centering
\begin{adjustbox}{max width=1\textwidth,center} 
\begin{tabular}{ccccccccccccc}
\hline
& & \multicolumn{5}{c}{$n=50$} && \multicolumn{5}{c}{$n=100$}
\\
& & \hhline{----- ~ -----}
 {$K$} &
{Approaches} & \multicolumn{1}{c}{$CvM$} & \multicolumn{1}{c}{$Dia$} & \multicolumn{1}{c}{$Cf_1$} & \multicolumn{1}{c}{$Cf_2$} & \multicolumn{1}{c}{$Kcop$} & &  \multicolumn{1}{c}{$CvM$} & \multicolumn{1}{c}{$Dia$} & \multicolumn{1}{c}{$Cf_1$} & \multicolumn{1}{c}{$Cf_2$} & \multicolumn{1}{c}{$Kcop$}
\\
& & \hhline{----- ~ -----}
%
%
\multirow{2}{*}{$2$} & $N(.2)$ & 4.7 & 5.2 & 6.2 & 6.1 & 6.0 && 4.7 & 3.9 & 4.7 & 4.7 & 4.0 \\  
 & $T(.2)$ & 4.1 & 3.9 & 4.8 & 4.6 & 6.0 & &4.4 & 4.6 & 5.2 & 5.3 & 5.0
 \\
& & \hhline{----- ~ -----}
%
%
\multirow{2}{*}{$3$} & $N(.2)$ & 3.3 & 4.8 & 5.7 & 4.8 & 4.0 && 2.9 & 4.1 & 3.3 & 3.8 & 4.0 \\  
 & $T(.2)$ & 3.0 & 4.4 & 3.9 & 3.7 & 6.0 && 4.0 & 5.1 & 4.3 & 4.5 & 6.0 \\
& & \hhline{----- ~ -----}
%
%
\multirow{2}{*}{$4$} & $N(.2)$ & 3.4 & 4.1 & 5.7 & 4.9 & 5.0 && 2.3 & 3.0 & 3.5 & 3.3 & 6.0 \\  
 & $T(.2)$ & 1.7 & 4.9 & 3.5 & 3.0 & 6.0 && 4.4 & 4.5 & 6.2 & 6.1 & 6.0 \\
\hline
%
\end{tabular}
\end{adjustbox}
\end{table}
\begin{table}
\caption{\label{TabQuessy} Empirical powers for  different alternatives studied in \cite{quessy16} with  $n=50$ and $n=100$.
}
    \centering
     \begin{adjustbox}{max width=1\textwidth,center}  
\begin{tabular}{ccccccccccccc}
\hline
%
&&  \multicolumn{5}{c}{$n=50$} & &\multicolumn{5}{c}{$n=100$}  \\
& & \hhline{----- ~ -----}
{$K$} & 
{Approaches} & \multicolumn{1}{c}{$CvM$} & \multicolumn{1}{c}{$Dia$} & \multicolumn{1}{c}{$Cf_1$} & \multicolumn{1}{c}{$Cf_2$} & \multicolumn{1}{c}{$Kcop$} & &  \multicolumn{1}{c}{$CvM$} & \multicolumn{1}{c}{$Dia$} & \multicolumn{1}{c}{$Cf_1$} & \multicolumn{1}{c}{$Cf_2$} & \multicolumn{1}{c}{$Kcop$} \\
%
%
& & \hhline{----- ~ -----}
\multirow{4}{*}{$2$} & $N(.4)$ & 16.8 & 15.5 & 21.3 & 19.8 & $\bm{32.0}$ & & 22.5 & 21.8 & 28.5 & 26.4 & $\bm{43.0}$ \\  
 & $T(.4)$ & 44.2 & 48.1 & $\bm{48.9}$ & 48.1 & 35.0 && 78.0 & 75.4 & $\bm{82.2}$ & 81.2 & 40.0 \\
 & $N(.6)$ & 51.3 & 48.9 & 62.2 & 58.1 & $\bm{68.0}$ & & 84.1 & 81.7 & 90.1 & 87.7 & $\bm{94.0}$ \\
 & $T(.6)$ & 88.6 & 92.2 & $\bm{90.8}$ & 90.7 & 60.0 & & $\bm{99.9}$ & $\bm{99.9}$ & $\bm{99.9}$ & $\bm{99.9}$ & 84.0 \\
& & \hhline{----- ~ -----}
%
%
\multirow{4}{*}{$3$} & $N(.4)$ & 12.0 & 12.7 & 17.8 & 16.1 & $\bm{64.0}$ & & 20.3 & 21.1 & 25.2 & 23.6 & $\bm{74.0}$ \\  
 & $T(.4)$ & 40.7 & 47.7 & 47.3 & 45.7 & $\bm{70.0}$ && 76.7 & 76.6 & $\bm{81.2}$ & 79.4 & 73.0 \\
 & $N(.6)$ & 47.3 & 48.2 & 63.0 & 56.8 & $\bm{92.0}$ && 85.8 & 84.7 & 91.5 & 88.6 & $\bm{99.0}$ \\
 & $T(.6)$ & 90.1 & $\bm{93.5}$ & 93.2 & 92.4 & 89.0 && 100.0 & 99.8 & 99.9 & $\bm{100.0}$ & 97.0 \\
& & \hhline{----- ~ -----}
\multirow{4}{*}{$4$} & $N(.4)$ & 9.8 & 12.0 & 16.1 & 14.1 & $\bm{78.0}$ & &19.6 & 20.6 & 27.1 & 25.0 & $\bm{84.0}$ \\  
 & $T(.4)$ & 34.6 & 41.7 & 43.9 & 42.3 & $\bm{84.0}$ && 74.4 & 75.5 & 78.8 & 77.6 & $\bm{86.0}$ \\
 & $N(.6)$ & 43.4 & 45.8 & 57.2 & 53.5 & $\bm{95.0}$ & & 81.5 & 80.3 & 88.9 & 86.4 & $\bm{100.0}$ \\
 & $T(.6)$ & 86.3 & 91.7 & 91.0 & 90.5 & $\bm{96.0}$ && 99.8 & 99.8 & $\bm{100.0}$ & $\bm{100.0}$ & 99.0 \\
\hline 
\end{tabular}
\end{adjustbox}
\end{table}

\FloatBarrier

\section{Conclusion}
In this paper, we introduced characteristic sequences, referred to as copula coefficients, for testing the equality of copulas. We developed a data-driven procedure in the two-sample case, accommodating both independent and paired populations. The extension to the $K$-sample case involves a second data-driven method, resulting in a two-step automatic comparison method. Our approach is applicable to all continuous random vectors, even in cases where the copula density does not exist.

Our method differs from the two-sample test proposed by \cite{remillard} and complements the $K$-sample test developed by \cite{quessy16, quessy21}, enabling the comparison of separate samples. The simulation study demonstrates the effectiveness of our approach, even for more than two populations. The test is user-friendly and performs efficiently.  We have limited our simulations to the case of ten samples, but larger dimensions are conceivable with this method. For future exploration, studying high dimensions within limited computation time may require dimension reduction by selecting a limited number of copula coefficients and vector components, which extends beyond the scope of this paper.

Comparing our method to existing approaches in the two-sample case, it appears as efficient as the competitor proposed by \cite{remillard}. In the $K$-sample case with $K>2$, numerical results suggest performance at least as good as those obtained by \cite{quessy16,quessy21}.  
In both cases of comparison, we used the previous models proposed by the authors. An R package of our procedure, named "Kcop," is available on CRAN.

Following the seminal work of \cite{quessy16} we can  adapt our procedure to test radial symmetry or exchangeability with a very similar statistic.
This idea is already nicely developed in \cite{quessy16,quessy2017, quessy21} with a general approach.

Eventually, our approach can be extended in various directions. Two potential directions include:
\bit
\item Copula coefficients can be used to obtain a simplified and unified  expression for some measures of association.
  Let us recall that for any continuous $d$-dimensional random variable $\bm{X}=(X_1,\ldots,X_d)$ with copula $C$,  one of the well-known popular multivariate versions of  Spearman's rho $\di\bm{\rho_{\bm{X}}}(C)$ can be expressed as 
(see \cite{Nelsen}): 
\be*
\di\bm{\rho_{\bm{X}}}(C)=h_{\rho}(d).\left\lbrace 2^{d}\int_{[0,1]^d}\pi(\bm{u}) dC(\bm{u})-1 \right\rbrace \text{ with } \di \pi(\bm u)=\prod_{j=1}^{d}u_j
\e*
and where $\di h_{\rho}(d)=\frac{d+1}{2^{d}-(d+1)}$. 
Then Spearman's rho coincides with the first copula coefficients, that is   
\be*
\di\bm{\rho_{\bm{X}}}(C)= h_{\rho}(d)\sum_{\m j \in \{0,1\}^{d}/\{0\}^{d}} \rho_{\m j}\prod_{k=1}^{d}\Big(\delta_{0,j_k}+\frac{\sqrt{3}}{3}\delta_{1,j_k}\Big).
\e*
For instance, for $d=3$, we have 
\be*
 \di\bm{\rho_{\bm{X}}}(C)= \frac{3}{4}\rho_{110} +\frac{3}{4}\rho_{101} +\frac{3}{4}\rho_{011}+\frac{3\sqrt{3}}{8}\rho_{111},
\e*
and we  deduce a novel estimator of the multivariate Spearman's rho as follows:
\be*
 \di \di\bm{\hat \rho_{\bm{X}}}(C) =
\di h_{\rho}(d)\sum_{\m j \in \{0,1\}^{d}/\{0\}^{d}} \w \rho_{\m j}\prod_{k=1}^{d}\Big(\delta_{0,j_k}+\frac{\sqrt{3}}{3}\delta_{1,j_k}\Big).
\e*

This estimator opens up possibilities for constructing tests comparing Spearman's rho. However, this requires the calculation of the asymptotic distributions of copula coefficients 
as proposed in  \cite{sinha1977multivariate}.

\item
Secondly, since the copula coefficients  characterize the dependence structure, we could  use such coefficients  for testing independence between random vectors in the same spirit as the penalized smooth tests proposed here. 
  \eit

\section*{Acknowledgements}

The authors would like to express their gratitude for the thorough reading, thoughtful comments, and numerous helpful suggestions provided by two anonymous referees and an Associate Editor. Their contributions greatly contributed to the improvement of this paper. The authors would like to extend special thanks to the Associate Editor for his helpful remarks, which led to Proposition 6. The second author would also like to acknowledge the support received from the Research Chair DIALog under the aegis of the Risk Foundation, an
initiative by CNP Assurances.




\bibliographystyle{imsart-number} 

\bibliography{LastrevisedEJS}

\appendix
\section{Proofs}
\label{appendproofs}
We detail the proof in the independent case. 
 The dependent case with $n_1=\cdots=n_K:=n$ is similar and will be indicated briefly in Appendix \ref{appendIndep}.  
Throughout the proofs, we  used the equality
$L_1(x)=\sqrt{3}(2x-1)$ and
the following inequalities are satisfied by Legendre polynomials (see \cite{abramowitz1964handbook}):
  \beg
  L_j(x) & \leq & c j^{1/2},  \quad \forall x \in [0,1]
  \label{ineg1}
  \\
  L'_j(x) & \leq & c' j^{5/2}, \quad \forall x \in [0,1],
   \label{ineg2}
 \en
   where $c>0$ and  $c'>0$ 
   are constant.
\subsection*{Proof of Proposition  \ref{prop00}}

%
%
%

From Corollary 6.7 of  \cite{moment2}, if $\mu $  is a Radon  measure  on $\R^p$  for which all moments
are finite and if there exists $\eps >0$ such that 
\begin{align}
\label{boundmoment}
    \int_{\R^p}e^{\eps \|x\|}\mu(dx) & < +\infty,
\end{align} 
then $\mu$ is said \textit{determinate}, that is:   if $\nu$ is a Radon measure with the same moments then $\nu = \mu$. 
Since $U$ is bounded on $[0,1]^p$, all its moments are finite and  (\ref{boundmoment}) is satisfied for all $\eps>0$. It follows that its distribution is {\it determinate}.   

\cqfd

\subsection*{Proof of Theorem \ref{thm1_ind}}
\label{appendproof2}
We want to show that $\P_0(D(\m n) >1) \rightarrow 0$ as $\m n$ tends to infinity.
We have
\begin{align}
 \nonumber
 \P_{0}\Big(D(\m n)> 1\Big)
 \nonumber
 &  \nonumber
 = \P_{0} \Big( \exists k \in \{2, \ldots,  d(\m n)\}:
  \nonumber
{V^{(1, 2)}_k -k\; p_{\m n} }\geq V^{(1, 2)}_1 -p_{\m n}\Big) \\
   & 
   =  \P_{0} \Big(  \exists k \in \{2, \ldots,  d(\m n)\}:
    \nonumber
 V_k^{(1, 2)} -V_1^{(1, 2)}\geq (k-1) p_{\m n}\Big) \\
  \nonumber
 &  
 =
 \P_{0} \Big(\exists k \in \{2, \ldots,  d(\m n)\}:
  \di \frac{n_1n_2}{n_1+n_2} \di\sum_{\m j \in {\cal H}^*(k)}    (r^{(1, 2)}_{\m j})^{2} \geq (k-1) p_{\m n}\big)
\nonumber
\\
 & 
 \leq
  \P_{0} \Big( \di \frac{n_1n_2}{n_1+n_2}   \di\sum_{\m j \in {\cal H}^*(d(\m n))}    (r^{(1, 2)}_{\m j})^{2} \geq {p_{\m n}} \Big),
 \label{ineg10}
 \end{align}
 with ${\cal H}(k)$ satisfying  (\ref{sumV}) and where  ${\cal H}^*(k)={\cal H}(k)\backslash {\cal H}(1)$. The last inequality comes from
 the fact that if a sum of $(k-1)$ positive terms, say
$\sum_{j=2}^k r_j $
is greater than a constant $c$, then necessarily there exists
a term $r_j$ such that $r_j > c/(k - 1)$.
The important point here is that $card({\cal H}^*(k))=k-1$, which  corresponds to the number of elements of the form  $(r_{\m j}^{(1,2)})^2 $ in the difference $V_k^{(1, 2)} -V_1^{(1, 2)}$.
For simplification of notation,  we  write ${\cal H}^*$ instead of ${\cal H}^*(d(\m n)) $.

Under the null $\rho^{(1)}_{\m j}= \rho^{(2)}_{\m j}$ and we  decompose $(r_{\m j}^{(1,2)})^2$ as follows
\beg
\label{egalr1}
\big(r_{\m j}^{(1, 2)}\big)^2 &=&
\big((\w \rho^{(1)}_{\m j}-\rho^{(1)}_{\m j}) - (\w {\rho}^{(2)}_{\m j} -\rho^{(2)}_{\m j})\big)^2
\\
\label{egalr2}
& \leq  &
2(\w \rho^{(1)}_{\m j}-\rho^{(1)}_{\m j})^2 + 2(\w { \rho}^{(2)}_{\m j} -\rho^{(2)}_{\m j})^2,
\en
that we combine with the
standard inequality for positive random variables: $\P(X+Y>z) \leq \P(X>z/2) + \P(Y>z/2)$, to get
\begin{align}
\label{AetB1}
\P_{0}\big(D(\m n)> 1\big) 
&\leq
\P_0 \big(    \di \frac{n_1n_2}{n_1+n_2} \sum_{\m j\in {\cal H}^*}(\w \rho^{(1)}_{\m j}-\rho^{(1)}_{\m j})^2 \geq  p_{\m n}/4\big)
  \nonumber
  \\
  &   \ \ \ 
  + \P_0\big(    \di \frac{n_1n_2}{n_1+n_2} \sum_{ \m j\in {\cal H}^*}(\w { \rho}^{(2)}_{\m j}-\rho^{(2)}_{\m j})^2 \geq  p_{\m n}/4\big)
  \\
  \label{AetB2}
  & :=  A + B.
  \end{align}
 We now study the first quantity $A$, the quantity $B$  being similar.
Writing
\be*
\tilde \rho^{(1)}_{\m j} & = &  \di\frac{1}{n_1} \di\sum_{s=1}^{n_1} L_{j_1}( U^{(1)}_{s,1})\cdots L_{j_p}( U^{(1)}_{s,p}),
\e*
we obtain
\beg
\label{decompose2}
\w \rho^{(1)}_{\m j}-\rho^{(1)}_{\m j}
& = &
(\w \rho^{(1)}_{\m j}-\tilde \rho^{(1)}_{\m j}) + (\tilde \rho^{(1)}_{\m j}-\rho^{(1)}_{\m j})
\ \  := \  \
E_{\m j} 
+G_{\m j},
\label{EetG}
\en
where
\begin{align*}
E_{\m j} & =   \di\frac{1}{n_1} \di\sum_{s=1}^{n_1}\left(  L_{j_1}( \w U^{(1)}_{s,1})\cdots L_{j_p}(\w  U^{(1)}_{s,p})  -L_{j_1}( U^{(1)}_{s,1})\cdots L_{j_p}( U^{(1)}_{s,p})
    \right) 
    \\
    G_{\m j} & =   \di\frac{1}{n_1} \di\sum_{s=1}^{n_1}\left(  L_{j_1}(  U^{(1)}_{s,1})\cdots L_{j_p}(  U^{(1)}_{s,p})  -\E\big( L_{j_1}( U^{(1)}_{1})\cdots L_{j_p}( U^{(1)}_{p}) \big)
    \right).
\end{align*}
Then we have
\begin{align}
\label{inegA}
A
& \leq 
\P_0\big(  \di \frac{n_1n_2}{n_1+n_2}   \sum_{ \m j\in {\cal H}^*}(E_{\m j})^2 \geq  p_{\m n}/16\big)
+\P_0\big(   \di \frac{n_1n_2}{n_1+n_2} \sum_{\m j\in {\cal H}^*}(G_{\m j})^2 \geq  p_{\m n}/16\big).
\end{align}
We first study the quantity involving $E_{\m j}$ in (\ref{inegA}).
Write
\beg
\label{supF}
S_{\ell}^{(1)} & = & \sup_x|\w F_{\ell}^{(1)}(x) - F_{ \ell}^{(1)}(x) |, \ \ \  \ell=1,\ldots, p.
\en
Applying the mean value theorem to $E_{\m j}$ we obtain
 \be*
| E_{\m j}| & \leq & \di\frac{1}{n_1} \di\sum_{s=1}^{n_1} \di\sum_{i=1}^p
S_{ i}^{(1)} \sup_x|L'_{j_i}(x)\di\prod_{u\neq i}L_{j_u}(x)|.
\e*
From (\ref{ineg1}) and (\ref{ineg2}) 
there exists a constant  $\t c>0$ such that
\beg
\label{inegal10}
|E_{\m j}|
& \leq  &
\t c \di \sum_{i=1}^p  S_{i}^{(1)} ({j_i}^{5/2}\di\prod_{u\neq i}{j_u}^{1/2}).
\en
When $\m j$ belongs to  ${\cal H}^* = {\cal H}^*(d(\m n))$ we  necessarily have $\| \m j \|_1 \leq d(n)$. Moreover $card({\cal H}^*)=d(\m n)-1$. It follows that
\beg
&&\P_0\big(    \di \frac{n_1n_2}{n_1+n_2} \sum_{\m j\in {\cal H}^*}(E_{\m j})^2 \geq   p_{\m n}/16\big)
\nonumber
\\
&&  \leq 
\P_0\big(    \di \frac{n_1n_2}{n_1+n_2} \sum_{\m j\in {\cal H}^*}\t c^2 \sum_{i=1}^p\sum_{i'=1}^p S_{ i}^{(1)}S_{ i'}^{(1)} j_{i}^{5/2}j_{i'}^{5/2} \di\prod_{s\neq i}{j_s}^{1/2} \di\prod_{s'\neq i'}{j_{s'}}^{1/2} \geq  p_{\m n}/16\big)
\nonumber
\\
&& 
  \leq 
\P_0\big(    \t c^2 \sum_{i=1}^p\sum_{i'=1}^p \di \frac{n_1n_2}{n_1+n_2} S_{ i}^{(1)}S_{i'}^{(1)}\sum_{\m j\in {\cal H}^*} d(\m n)^{p+4} \geq  p_{\m n}/16\big)
\nonumber
\\
&& 
  \leq 
\P_0\big(    \t c^2 \sum_{i=1}^p\sum_{i'=1}^p \di \frac{n_1n_2}{n_1+n_2} S_{ i}^{(1)}S_{i'}^{(1)} d(\m n)^{p+5} \geq  p_{\m n}/16\big)
\nonumber
\\
&&   \rightarrow 0 \ {\rm \ as  } \ \m n \rightarrow \infty,
\label{conv3}
\en
 since for all $i=1,\ldots, p$, $\sqrt{\di \frac{n_1n_2}{n_1+n_2}}S_{i}^{(1)}$ converges in law to a Kolmogorov distribution and  $d(\m n)^{p+5} = o(p_{\m n})$ by {\bf(A)}.

Coming back to (\ref{EetG}), we now study the quantity involving $G_{\m j}$. First note that $\E(G_{\m j})=0$.
Moreover,
$\V(G_{\m j}) = \E((G_{\m j})^2) = \V(\di\prod_{i=1}^pL_{j_i}(U^{(1)}_i))/n_{1}$. Then, by Markov inequality we have
\be*
\P_0\big(    \di \frac{n_1n_2}{n_1+n_2} \sum_{\m j\in {\cal H}^*}(G_{\m j})^2 \geq  p_{\m n}/16\big)
 &\leq&
\di \frac{n_2}{n_1+n_2}\di\frac{ \di   \sum_{\m j\in {\cal H}^*}\V(\di\prod_{i=1}^p L_{j_i}(U^{(1)}_i))} { p_{\m n}/16}
\e*
and from (\ref{ineg1})
there exists a constant $c>0$ such that
\be*
\V\big(\di\prod_{i=1}^p L_{j_i}(U^{(1)}_i)\big) &\leq&
c^2 \di\prod_{i=1}^p j_i.
\e*
It follows that
\begin{align}
\P_0\big(    \di \frac{n_1n_2}{n_1+n_2} \sum_{\m j\in {\cal H}^*}(G_{\m j})^2 \geq  p_{\m n}/16\big)
\;
\leq
\di \frac{n_2}{n_1+n_2}\di\frac{c^2 d(\m n)^p} { p_{\m n}/16}
&&\rightarrow  0, \ {\rm as  } \ \m n \rightarrow \infty.
\label{conv4}
\end{align}
We now  combine (\ref{conv3}) and (\ref{conv4}) with (\ref{inegA}) to conclude that
$ A \rightarrow 0$,   as    $\m n \rightarrow \infty$.

In the same manner we can show that
$
 B \rightarrow 0, \ {\rm as  } \ \m n \rightarrow \infty$,
which completes the proof.

 \cqfd


   \subsection*{Proof of Theorem  \ref{thm2_ind}}

Let $\m j=(1,1,\ldots, 0,0)$. We have
$
 \di V^{(1, 2)}_1 = T^{(1, 2)}_{2,1}= \left(\sqrt{\di \frac{n_1n_2}{n_1+n_2}}r^{(1,2)}_{\m j}\right)^2
$
and we  can decompose $\sqrt{\di \frac{n_1n_2}{n_1+n_2}}r^{(1,2)}_{\m j}$ under the null as follows:
\begin{eqnarray}
\nonumber
\sqrt{\di \frac{n_1n_2}{n_1+n_2}}r^{(1,2)}_{\m j}&=&\sqrt{\di \frac{n_1n_2}{n_1+n_2}}\left(\w \rho^{(1)}_{\m j}-\w \rho^{(2)}_{\m j}\right) 
\\
\nonumber
 &=&
\sqrt{\di \frac{n_1n_2}{n_1+n_2}}\Bigg( \di\frac{1}{{n_1}}\di\sum_{i=1}^{n_1}
L_{1}(\w U^{(1)}_{i, 1})L_{1}(\w U^{(1)}_{i,2})\\\nonumber
&&\qquad-
\di\frac{1}{{n_2}}\di\sum_{i=1}^{n_2}
L_{1}(\w U^{(2)}_{i, 1})L_{1}(\w U^{(2)}_{i,2})\Bigg)\\
\label{decomp1}
& =&
\sqrt{\di \frac{n_1n_2}{n_1+n_2}}\left( \di\frac{1}{{n_1}}(\di\sum_{i=1}^{n_1}
L_{1}(\w U^{(1)}_{i, 1})L_{1}(\w U^{(1)}_{i,2}) -m\right)
\nonumber
\\
&& \ \ \ \ 
-\sqrt{\di \frac{n_1n_2}{n_1+n_2}}\left(\di\frac{1}{{n_2}}(\di\sum_{i=1}^{n_2}
L_{1}(\w U^{(2)}_{i, 1})L_{1}(\w U^{(2)}_{i,2})-m\right)
\\
\label{decomp2}
& : = &
R^{(1)}_{\m n} - R^{(2)}_{\m n},
\end{eqnarray}
where, under the null
\begin{align*}
 m =\E(L_{1}( U^{(1)}_{i, 1})L_{1}( U^{(1)}_{i,2})) =\E(L_{1}( U^{(2)}_{i, 1})L_{1}( U^{(2)}_{i,2})).
\end{align*}
By Taylor expansion, using the fact that the Legendre polynomials satisfy  $L_1'=2\sqrt{3}$ and $L_1''=0$,  we obtain
\be*
R_{\m n}^{(1)}
& = &
\sqrt{\di \frac{n_1n_2}{n_1+n_2}} 	\Bigg\lbrace
\sum_{i=1}^{n_1}
\big(L_{1}( F_{1}^{(1)}(x_{i,1}^{(1)}))L_{1}( F^{(1)}_{2}(X_{i,2}^{(1)}))
d\w {\m F}^{(1)}(x,y) - m \big)
\\
& & \  \ +
\di\int\di\int
(\w F_{1}^{(1)}(x)-F_{1}^{(1)}(x))2\sqrt{3} L_{1}( F^{(1)}_{2}(y))
d \m F^{(1)}(x,y)
\\
&& \ \  +
\di\int\di\int
(\w F_{2}^{(1)}(y)-F_2^{(1)}(y)) 2\sqrt{3}L_{1}( F^{(1)}_{1}(x))
d \m F^{(1)}(x,y)
\\
&&  
+
\di\int\di\int
(\w F_{1}^{(1)}(x)-F_{1}^{(1)}(x)) 2\sqrt{3}L_{1}( F^{(1)}_{2}(y))
d\big(\w {\m F}^{(1)} - \m F^{(1)}\big)(x,y)
\\
&& +
\di\int\di\int
(\w F_{2}^{(1)}(y)-F_2^{(1)}(y)) 2\sqrt{3}L_{1}( F^{(1)}_{1}(x))
d\big(\w {\m F}^{(1)} - \m F^{(1)}\big(x,y)\Bigg\rbrace
\\
&:= &
\sqrt{\di \frac{n_1n_2}{n_1+n_2}} \left( A^{(1)}_{1, n_1} + A^{(1)}_{2, n_1} + A^{(1)}_{3, n_1} + B^{(1)}_{ n_1} + C^{(1)}_{ n_1} \right).
\e*
By symmetry,  the second term $R^{(2)}_{\m n}$ can be expressed as:
\begin{align*}
    R^{(2)}_{\m n}
    &=
    \sqrt{\di \frac{n_1n_2}{n_1+n_2}} \left( A^{(2)}_{1, n_2} + A^{(2)}_{2, n_2} + A^{(2)}_{3, n_2} + B^{(2)}_{ n_2} + C^{(2)}_{ n_2} \right)
\end{align*}
and finally
\be*
    \sqrt{\di \frac{n_1n_2}{n_1+n_2}}\; r^{(1,2)}_{\m j}
    & = &
     \sqrt{\di \frac{n_1n_2}{n_1+n_2}} \left( A^{(1)}_{1, n_1} + A^{(1)}_{2, n_1} + A^{(1)}_{3, n_1}-A^{(2)}_{1, n_2} \right. 
     \\
     & & \hspace*{1.1cm}\left. - A^{(2)}_{2, n_2} - A^{(2)}_{3, n_2} + B^{(1)}_{ n_1} + C^{(1)}_{ n_1} - B^{(2)}_{ n_2} - C^{(2)}_{ n_2} \right).
\e*
%
%
%
This expression is very similar to the expansion used in \cite{sinha1977multivariate} (see his proof of Theorem 1) and \cite{bhuchongkul1964class} (see his equation (3.4)). We imitate their approach here. Therefore, we will show that
$\di \sqrt{n_1}\sum_{i=1}^{3} A^{(1)}_{i, n_1} $ and $\di \sqrt{n_2}\sum_{i=1}^{3} A^{(2)}_{i, n_2} $ have a limiting normal distribution and the rest of the terms are all $o_{\P}(1)$.  
Using the expression of the empirical cdf we can rewrite
\be*
A^{(1)}_{1,\m n} + A^{(1)}_{2,\m n} &+& A^{(1)}_{3,\m n}
= 
\di\frac{1}{n_1}\di\sum_{i=1}^{n_1}    {\large \Big\{ }
L_{1}(F_{1}^{(1)}(X^{(1)}_{1,i})L_{1}(F^{(1)}_{2}(X^{(1)}_{2,i}))
-m 
\\
&&
+ 2\sqrt{3}\di\int\di\int (\mathds{1}(X^{(1)}_{1,i}\leq x )-F_{1}^{(1)}(x)) L_{1}( F^{(1)}_{2}(y))
d \m F^{(1)}(x,y)
\\
&&
\left.
+
2\sqrt{3}\di\int\di\int (\mathds{1}(X^{(1)}_{2,i}\leq y )-F_{2}^{(1)}(y)) L_{1}( F^{(1)}_{1}(x))
d \m F^{(1)}(x,y)
\right\}
\\
 &:= &
\di\frac{1}{n_1}\di\sum_{i=1}^{n_1}  (Z^{(1)}_{1,i} + Z^{(1)}_{2,i} + Z^{(1)}_{3,i}) \ : = \di\frac{1}{n_1}\di\sum_{i=1}^{n_1}  Z^{(1)}_{i},  
\e*
%
%
%
where $Z^{(1)}_i$ are iid random variables. 
By symmetry we get 
\be*
A^{(2)}_{1,\m n} + A^{(2)}_{2,\m n} + A^{(2)}_{3,\m n}
& = &
\di\frac{1}{n_2}\di\sum_{i=1}^{n_2}  (Z^{(2)}_{1,i} + Z^{(2)}_{2,i} + Z^{(2)}_{3,i}) \ : = \di\frac{1}{n_2}\di\sum_{i=1}^{n_2}  Z^{(2)}_{i}.  
\e*
Clearly $\E(Z^{(1)}_{1,i}-Z^{(2)}_{1,i}) = 0$.
Since  $ \E(\mathds{1}(X^{(1)}_{1,i}\leq x ))=F_{1}^{(1)}(x)$ and $ \E(\mathds{1}(X^{(2)}_{1,i}\leq x ))=F_{1}^{(2)}(x)$,
we also have  $\E(Z^{(1)}_{2,i}-Z^{(2)}_{2,i})=0$ and similarly   $\E(Z^{(1)}_{3,i}- Z^{(2)}_{3,i})= 0$.
Moreover, $Z_i^{(1)}$ and $Z_i^{(2)}$ have finite variances. 
Applying the Central Limit Theorem to the independent iid series $Z_i^{(1)}$ and $Z_i^{(2)}$ we obtain
\be*
\sqrt{\di \frac{n_1n_2}{n_1+n_2}}\left(
\di\frac{1}{n_1}\di\sum_{i=1}^{n_1}  Z^{(1)}_{i}+ \di\frac{1}{n_2}\di\sum_{i=1}^{n_2}  Z^{(2)}_{i}
\right)
& \rightarrow & N(0, \si^2(1,2)),
\e*
with
\be*
\si^2(1,2) & = &
(1-a_{1,2})\V(Z_i^{(1)}) 
+a_{1,2}\V(Z_i^{(2)} )
\e*
where $a_{1,2}$ is given by (\ref{ratiorate}), 
    and where, for $j=1,2$,  
\be* 
\V(Z_i^{(j)})  &=& \V
\Bigg(
L_1(U^{(j)}_{1}) L_1(U^{(j)}_{2})
\\
&&
\ \ \ \ + 2\sqrt{3}\di\int\di\int \big(\mathds{1}(X^{(j)}_{1}\leq x )-F_{1}^{(j)}(x)\big)L_{1}( F^{(j)}_{2}(y))
d \m F^{(j)}(x,y)
\\
&&
\  \  \ \  +2\sqrt{3}\di\int\di\int \big(\mathds{1}(X^{(j)}_{2}\leq y )-F_{2}^{(j)}(y)\big) L_{1}( F^{(j)}_{1}(x))
d \m F^{(j)}(x,y)
\Bigg).
\e*
We now proceed to check  that $B^{(s)}_{ n_s}$, $C^{(s)}_{n_s}$ are $o_{\P}(n_s^{-1/2})$ for $s=1,2$. 
The asymptotic negligibility of $B^{(s)}_{ n_s}, s=1,2$ and $C^{(s)}_{n_s},s=1,2$ follows directly from  those of $\bm{B_{1N}}$ and $\bm{B_{2N}}$ in \cite{bhuchongkul1964class}. The arguments are exactly similar
to those of \cite{bhuchongkul1964class} (see his proof of Theorem 1) and we therefore omit the details.


 \cqfd

  \subsection*{Proof of Proposition \ref{prop1_ind}}

Let us define
\be*
\di \overline{W}^{(s)}=\frac{1}{n_s}\sum_{i=1}^{n_s}W_{i}^{(s)}, \ \ \text{for }s=1,2,
\e*
 where
\be*
W_{i}^{(s)} &=& L_1(U^{(s)}_{i,1}) L_1(U^{(s)}_{i,2})\\
 &&\quad\quad + 2\sqrt{3}\di\int\di\int \big(\mathds{1}(X^{(s)}_{i,1}\leq x )-F_{1}^{(s)}(x)\big)L_{1}( F^{(s)}_{2}(y))
d F^{(1)}(x,y)
\\
 &&\quad\quad
+2\sqrt{3}\di\int\di\int \big(\mathds{1}(X^{(s)}_{i,2}\leq y )-F_{2}^{(s)}(y)\big) L_{1}( F^{(s)}_{1}(x))
d F^{(1)}(x,y).
\e*
%
We focus on $s=1$, the case $s=2$ being similar. We have
\beg\label{prob1 2}
\di \frac{1}{n_1}\di\sum_{i=1}^{n_1}\Big(W_{i}^{(1)} - \overline{W}^{(1)}\Big)^2 \overset{\P}{\longrightarrow} \sigma^2{(1)}.
\en
According to  Slusky's Lemma and (\ref{prob1 2}), the proof is completed by showing that
\be*
\di \frac{1}{n_1}\di\sum_{i=1}^{n_1}\Big(W_{i}^{(1)}- \overline{W}^{(1)}\Big)^2 -\w\sigma^2{(1)}\overset{\P}{\longrightarrow} 0.
\e*
We have
\be*
&&\di \frac{1}{n_1}\di\sum_{i=1}^{n_1}\Big(W_{i}^{(1)} - \overline{W}^{(1)} \Big)^2 -\w\sigma^2{(1)}\\
&\quad=& \di \frac{1}{n_1}\di\sum_{i=1}^{n}\Big(W_{i}^{(1)} \Big)^2-\Big(\overline{W}^{(1)}\Big)^2 -
\di \frac{1}{n_1}\di\sum_{i=1}^{n_1}\Big(M_{i}^{(1)}\Big)^2 +\Big(\overline{M}^{(1)}\Big)^2\\
&=&
\di \frac{1}{n_1}\sum_{i=1}^{n_1}\Big(W_{i}^{(1)} -M_{i}^{(1)} \Big)\Big(W_{i}^{(1)} +M_{i}^{(1)} -\overline{M}^{(1)}-\overline{W}^{(1)}\Big).
\e*
From (\ref{ineg1}), there exists a constant $\kappa>0$ such that, for all $n_1>0$ and for all $i=1,\ldots, n_1$,
\be*
\di \max(|W_{i}^{(1)}|,|M_{i}^{(1)}|)\leq \kappa,
\e*
which implies  that
\be*
\nonumber
\di\left|  \frac{1}{n_1}\di\sum_{i=1}^{n_1}\Big(W_{i}^{(1)} - \overline{W}{^(1)}\Big)^2 -\w\sigma^2{(1)} \right|
\leq
 \frac{8\kappa}{n_1}\sum_{i=1}^{n_1}\left|W_{i}^{(1)}-M_{i}^{(1)}  \right| .
\e*
It remains to prove that  $\di W_{i}^{(1)}-M_{i}^{(1)}  \overset{\P}{\longrightarrow} 0$.
We have
$W_{i}^{(1)} - M_{i}^{(1)} = I_{i,1}  + 2\sqrt{3}I_{i,2} + 2\sqrt{3}I_{i,3}$,
%
where
\be*
I_{i,1} &=&
L_{1}( U^{(1)}_{i,1})L_{1}( U^{(1)}_{i,2})-L_{1}(\w U^{(1)}_{i,1})L_{1}(\w U^{(1)}_{i,2})
\\
I_{i,2} &=&
\int\left(\mathds{1}(X^{(1)}_{i,1}\leq x)-
F^{(1)}_{1}(x)\right)L_{1}\big( F^{(1)}_2(y)\big)dF^{(1)}(x,y)
\\
 &&\qquad\qquad\qquad\qquad\quad-\frac{1}{n_1}\di\sum_{k=1}^{n_1}\Big(\mathds{1}(X^{(1)}_{i,1}\leq X^{(1)}_{k,1})-
\w U^{(1)}_{k,1}\Big)L_{1}(\w U^{(1)}_{k,2})
\\
I_{i,3} &=&
\int\left(\mathds{1}(X^{(1)}_{i,2}\leq x)-
F^{(1)}_{1}(x)\right)L_{1}\big( F^{(1)}_2(y)\big)dF^{(1)}(x,y)
\\
 &&\qquad\qquad\qquad\qquad\quad-\frac{1}{n_1}\di\sum_{k=1}^{n_1}\Big(\mathds{1}(X^{(1)}_{i,2}\leq X^{(1)}_{k,2})-
\w U^{(1)}_{k,2}\Big)L_{1}(\w U^{(1)}_{k,1}).
\e*
%
%
%
%
Since $\di L_{1}(t)=\sqrt{3}(2t-1)$, we get
\be*
|\di I_{i,1} |&=&
\big|2\sqrt{3}L_{1}( U^{(1)}_{i,1})\Big( U^{(1)}_{i,2}-
\w U^{(1)}_{2}\Big)+
2\sqrt{3}L_{1}( \w U^{(1)}_{2})\Big( U^{(1)}_{i,1}-\w U^{(1)}_{1}\Big)\big|\\
&\leq & 6 (S^{(1)}_{2} + S^{(1)}_{1}) \ \
= o_{\P}(1),
\e*
where $S^{(1)}_{2}$ and $S^{(1)}_{1}$ are given by (\ref{supF}).
We next show that  $\di I_{i,2} = o_{\P}(1) $.  We have 
\be*
\di I_{i,2} &=&
\frac{1}{n_1}\di\sum_{i=1}^{n_1}U^{(1)}_{i,1}L_{1}( U^{(1)}_{i,2})
-\iint F^{(1)}_{1}(x)L_{1}\big( F^{(1)}_{2}(y)\big)dF^{(1)}_{1,2}(x,y)\\
&& \ +
 \iint\mathds{1}_{X^{(1)}_{k,1}\leq x^{(1)}_{1}}L_{1}\big( F^{(1)}_{2}(y)\big)dF^{(1)}_{1,2}(x,y)-\frac{1}{n_1}
\di\sum_{i=1}^{n_1}\mathds{1}_{X^{(1)}_{k,1}\leq X^{(1)}_{i,1}}L_{1}(U^{(1)}_{i,2})\\\
&:= & I_{2,k}^{1} +I_{2,k}^{2}.
\e*
To deal with $\di I_{2,k}^{1}$, we note that
\be*
I_{2,k}^{1}&=& 
\frac{1}{n_1}\di\sum_{s=1}^{n_1}U^{(1)}_{s,1}L_{1}( U^{(1)}_{s,2})
-\E\Big(U^{(1)}_{1}L_{1}( U^{(1)}_{2})\Big).
\e*
Since the random vectors $(\m U_{i}^{(1)} : = \di (U^{(1)}_{i,1},U^{(1)}_{i,2}))_{i=1,\ldots, n_1}$ 
are iid, 
the
weak law of large numbers  and the continuous mapping theorem show that
\be*\label{I i 221}
\di I_{2,k}^{1}=o_{\P}(1).
\e*
For  $\di I_{2,k}^{2}$, we can write 
\be*
\di I_{2,k}^{2}
&=&
 \iint\mathds{1}_{F^{(1)}_{1}(X^{(1)}_{k,1})\leq F^{(1)}_{1}(x)}L_{1}\big( F^{(1)}_{2}(y)\big)dF^{(1)}_{1,2}(x,y)\\
 &&\qquad\qquad\qquad\qquad\quad-
 \frac{1}{n_1}
\di\sum_{i=1}^{n_1}\mathds{1}_{F^{(1)}_{1}(X^{(1)}_{k,1})\leq F^{(1)}_{1}(X^{(1)}_{i,1})}L_{1}(U^{(1)}_{i,2})\\
&=&
\int_{0}^{1}\int_{0}^{1} \mathds{1}_{U^{(1)}_{k,1}\leq u}L_{1}(v)dC^{(1)}(u,v)
-\frac{1}{n_1}\di\sum_{i=1}^{n_1}\mathds{1}_{U^{(1)}_{k,1}\leq U^{(1)}_{i,1}}L_{1}(U^{(1)}_{i,2})
\e*
and since  $U^{(1)}_{i,1}$ has continuous uniform distribution it follows that
\be*
\di | I_{2,k}^{2}|&\leq &
\sup_{t\in [0,1]}\Bigg\lvert
\frac{1}{n_1}\di\sum_{i=1}^{n_1}\mathds{1}_{t\leq U^{(1)}_{i,1}}L_{1}(U^{(1)}_{i,2})\\
 &&\qquad\qquad\qquad\qquad\quad-\int_{0}^{1}\int_{0}^{1} \mathds{1}_{t\leq u^{(1)}_{1}}L_{1}(u^{(1)}_{2})dC^{(1)}(u^{(1)}_{1},u^{(1)}_{2})
 \Bigg\rvert\\
&\leq &
\sup_{t\in [0,1]}\left|
\frac{1}{n_1}\di\sum_{i=1}^{n_1}\mathds{1}_{t\leq U^{(1)}_{i,1}}L_{1}(U^{(1)}_{i,2})
-\E\left(\mathds{1}_{t\leq U^{(1)}_{1}}L_{1}(U^{(1)}_{2})\right)
 \right|\\
 &\leq &
\sup_{t\in [0,1]}\left|
g\left(t,\m U^{(1)}_{1},\ldots,\m U^{(1)}_{n_1}\right)
-\E\Big(g\big(t,\m U^{(1)}\big)\Big)
 \right|,
\e*
where
%
\be*
\di g\left(t,z_1,\ldots,z_{n_1}\right)=\frac{1}{n_1}
\di \sum_{k=1}^{n_1}\mathds{1}_{t\leq u_{k}}L_{1}(v_{k}), \text{and } z_k=(u_{k},v_{k})\text{ for }k=1,\ldots,n_1.
\e*
Observe that for all $\di t\in [0,1]$,
\be*
\sup_{\substack{z_{1},\ldots,z_{n_1},\\z'_{i}}}\left|g(t,z_{1},\ldots,z_{n_1})
-g(t,z_{1},\ldots,z_{i-1},z'_{i},z_{i+1},\ldots,z_{n_1})\right| \leq \frac{2\|L_{1}\|_{\infty}}{n_1},
\e*
with $\frac{2\|L_{1}\|_{\infty}}{n_1}=\frac{4\sqrt{3}}{n_1}$,
that is,  if we change the $i$th variable $z_i$ of $g$ while keeping all the others fixed, then the value of the
function does not change by more than $4\sqrt{3}/n_1$.
Then, by McDiarmid's inequality, we get $\forall \epsilon>0$
\be*
\P\left(\forall t,  \left|
g\left(t,\m U^{(1)}_{1},\ldots,\m U^{(1)}_{n_1} \right)
-\E\Big(g\big(t,\m U^{(1)} \big)\Big)
 \right|\geq \epsilon\right)
&  \leq &2e^{-n_1\epsilon^{2}/24} \underset{n_1\to \infty}{\longrightarrow} 0.
\e*
It implies that
%
$\di I_{2,k}^{2} =o_{\P}(1)$,
%
and we conclude that 
$\di I_{i,2}=o_{\P}(1)$ and  similarly that $\di I_{i,3}=o_{\P}(1)$.
It follows that
 $\di W_{i}^{(1)}-M_{i}^{(1)} \overset{\P}{\longrightarrow} 0$  
 which completes the proof.

\cqfd



 \subsection*{Proof of Theorem \ref{thm3_ind}}

Let us   prove  that
$\P(s( \m n) \geq 2)$ vanishes as $\m  n\rightarrow +\infty$.
 By definition of $s(\m n)$  we have:
\begin{align*}
\nonumber
\P\big(s(\m n)\geq 2\big)
&= \P\big( \mbox{there exists}~ 2\leq k\leq v(K):~V_k - k p_{\m n} \geq V_1 - p_{\m n} \big)
\nonumber \\
&= \P\big( \mbox{there exists}~ 2\leq k\leq v(K): V_k-V_1  \geq (k -1) p_{\m n} \big)
\nonumber
\\
&= \P\big( \exists~ 2\leq k\leq v(K): \di\sum_{2 \leq ord_{\cal V}(\ell,m)\leq k} V^{(\ell,m)}_{D(\m n)} \geq (k -1) p_{\m n} \big).
\end{align*}
Since the previous sum contains $(k-1)$ positive elements, there is at least one element greater than $p_{\m n}$. It follows that
\begin{align*}
\P\big(s(\m n)\geq 2\big)
&\leq \P \big( \exists (\ell,m)~ \mbox{with}~ {2 \leq ord_{\cal V}(\ell,m) \leq  v(K)}:~  V^{(\ell,m)}_{D(\m n)} \geq p_{\m n} \big) \nonumber
\\
&\leq \P\left (\di\sum_{2\leq  ord_{\cal V}(\ell,m)\leq v(K)} V^{(\ell,m)}_{D(\m n)}  \geq p_{\m n} \right ).
\end{align*}
 First, we can remark that ${\cal V}(K)$ is finite and  then there is a finite number of terms in
 $\di \sum_{2\leq  ord_{\cal V}(\ell,m)\leq v(K)}V^{(\ell,m)}_{D(\m n)}  $. It follows that we simply have to show that the probability
$\P( V^{(\ell,m)}_{D(\m n)}  \geq p_{\m n})$  vanishes as $ \m n\rightarrow +\infty$ for any values of $(\ell,m)$ .
 Since $D(\m n) \leq d(\m n)$ have:
\beg
\P( V^{(\ell,m)}_{D(\m n)}  \geq p_{\m n})
&\leq &
\P( V^{(\ell,m)}_{d(\m n)}  \geq p_{\m n})
\nonumber
\\
\label{egal1000}
 &=&
 \P_{0} \Big(  \di \frac{n_{\ell}n_m}{n_{\ell}+n_m}   \di\sum_{\m j \in {\cal H}(d(\m n))} (r^{(\ell, m)}_{\m j})^{2} \geq p_{\m n} \Big).
 \en
Comparing (\ref{egal1000})  and (\ref{ineg10})  we can see that the study  is now similar in spirit to the two-sample case  and
we can simply  mimic the proof of Theorem \ref{thm1_ind} to conclude.

 \cqfd

 \subsection*{Proof of Proposition \ref{prop3_ind}}
\label{appendproof1}
We give the proof for  the case $k>1$, the particular case $k=1$ being similar. For simplification of notation,  we  now write ${\cal H}$ instead of ${\cal H}(d(\m n)) $. We first show that $\P(s( \m n) \geq k)$ tends to 1. Under $H_1(k)$,   we have for all $k'<k$:
\begin{align}
\nonumber
\P(s(\m n) < k )
& \leq
\P\left (V_k - k p_{\m n} \leq V_{k'} - k'p_{\m n} \right )\\
\nonumber
&
= 1 -
\P\left ((V_k - V_{k'})  \geq  (k-k') p_{\m n} \right )
\\
\nonumber
&
=
1-\P\left (\di \sum_{ k'< r_{\cal V}(\ell,m) \leq k}  V^{(\ell,m)}_{D(\m n)}  \geq  (k-k') p_{\m n} \right )
\\
\nonumber
& = 1-
\P\left (\di \sum_{ k'< r_{\cal V}(\ell,m) \leq k} \di \frac{n_{\ell}n_m}{n_{\ell}+n_m}  \di\sum_{\m j \in {\cal H}}
(r_{\m j }^{(\ell,m)})^2
\geq  (k-k') p_{\m n} \right )
\\
\label{1moinsp}
&
\leq
1 -
\P\left (
\mathds{1}_{\left\{ r_{\cal V}(\ell,m)=k \right\}}\,
\di \frac{n_{\ell}n_m}{n_{\ell}+n_m}   \sum_{\m j \in {\cal H}}
(r_{\m j }^{(\ell,m)})^2
\geq  (k-k') p_{\m n} \right ).
\end{align}
 When $r_{\cal V}(\ell,m)=k$, under $H_1(k)$, since $C^{(\ell)}\neq C^{(m)}$,   there exists $\m j_0$ such that
$\rho^{(\ell)}_{\m j_0} \neq \rho^{(m)}_{\m j_0}$.
 We have
 \begin{align}
\P\left (
\mathds{1}_{\left\{ r_{\cal V}(\ell,m)=k \right\}}
\di \frac{n_{\ell}n_m}{n_{\ell}+n_m}   \sum_{\m j \in {\cal H}}
(r_{\m j }^{(\ell,m)})^2
\geq  (k-k') p_{\m n} \right )
\hspace*{1.5cm}& 
\nonumber
\\
\label{inequality10}
  \geq
\P\left (
\mathds{1}_{\left\{ r_{\cal V}(\ell,m)=k \right\}}
\di \frac{n_{\ell}n_m}{n_{\ell}+n_m}   \mathds{1}_{\m j_0 \in {\cal H}}
(r_{\m j_0 }^{(\ell,m)})^2
\geq  (k-k') p_{\m n} \right )&,
 \end{align}
and we can decompose $r_{\m j_0 }^{(\ell,m)}$ as follows
\begin{align}
r_{\m j_0 }^{(\ell,m)}
& =
\left( (\w \rho^{(\ell)}_{\m j_0}
-\rho^{(\ell)}_{\m j_0})
- (\w \rho^{(m)}_{\m j_0} -\rho^{(m)}_{\m j_0}) \right)
+\left(\rho^{(\ell)}_{\m j_0} -\rho^{(m)}_{\m j_0})\right)
\label{egalR00}
\ \ := \ \ (A-B)+D.
\end{align}
 We first decompose the quantities $A$ and $B$. We only detail the calculus for $A$, since the case of $B$ is similar.  We have
 \be*
\label{decompose21}
A
& = &
(\w \rho^{(\ell)}_{\m j_0}-\tilde \rho^{(\ell)}_{\m j_0}) + (\tilde \rho^{(\ell)}_{\m j_0}-\rho^{(\ell)}_{\m j_0})
\nonumber
\ \ :=  \ \
E_{\m j_0} 
+G_{\m j_0}.
\e*
 We can reuse (\ref{inegal10}) to get:
 \be*
|E_{\m j_0}|
& \leq  &
\t c \di \sum_{i=1}^p  S_{i}^{(\ell)} ({j_i}^{5/2}\di\prod_{u\neq i}{j_u}^{1/2})
\ \ \leq \ \
\t c' \|\m j_0\|_1^{(p+4)/2}\di \sum_{i=1}^p  S_{i}^{(\ell)},
\e*
for some constants $\t c $ and $\t c'$.
Since $\sqrt{n_{\ell}}S_{i}^{(\ell)} = O_{\P}(1)$ (see for instance \cite{massart1990tight}) we have
$     n_{\ell}
E_{\m j_0}^2  = O_{\P}( 1) $.
As $G_{\m j_0}$ is an empirical estimator we also have $     n_{\ell}
G_{\m j_0}^2  = O_{\P}( 1)$, which yields
\begin{align}
\label{quantityA}
    n_{\ell} A^2 & =  O_{\P}( 1).
\end{align}

We now consider the quantity $D$ in (\ref{egalR00}). The inequality  $\rho^{(\ell)}_{\m j_0} \neq \rho^{(m)}_{\m j_0}$ implies that
 \begin{align}
\label{quantityD} \di \frac{n_{\ell}n_m}{n_{\ell}+n_m}  D^2
  = O(\m n).
   \end{align}
Finally, under $H_1(k)$, we combine (\ref{quantityA}) and (\ref{quantityD}) with (\ref{egalR00}) to get\\
$ \di \frac{n_{\ell}n_m}{n_{\ell}+n_m} 
(r_{\m j_0 }^{(\ell,m)})^2
     =
    O_{\P}(\m n)$.
If we prove that $\mathds{1}_{\m j_0 \in {\cal H}(D(\m n))} \rightarrow 1$ as $ n$ tends to infinity
then   (\ref{inequality10}) tends to 1,   from assumption {\bf (B)}.
Mimicking the proof of Theorem \ref{thm1_ind} we can prove that
$\P(D(\m n)< ord(\m j_0,\|\m j_0\|_1)) \rightarrow 0 $ which gives the result.

Our next goal is to determine the limit of $\P(V < \eps)$ for $\eps >0$.
It is sufficient to prove that $\P(V_{s(\m n)} < \eps) \rightarrow 0$ as $\m n$ tends to infinity.
We have
\begin{align*}
    \P(V_{s(\m n)} < \eps)
    & =
    \di\sum_{s=1}^{v(K)}
    \P(V_{s} < \eps \cap
    s(\m n) = s)
    \\
    & =
      \di\sum_{s=1}^{k-1}
  \P(V_{s} < \eps \cap
    s(\m n) = s)
    +
     \di\sum_{s=k}^{v(K)}
  \P(V_{s} < \eps \cap
    s(\m n) = s)
    \\
    & \leq
       \di\sum_{s=1}^{k-1}
  \P(V_{s} < \eps \cap
    s(\m n) = s)
    +
     \di\sum_{s=k}^{v(K)}
  \P(V_{s} < \eps )
  \ \  := E +F .
\end{align*}
From what has already been proved, under $H_1(k)$
\begin{align*}
\lim_{\m n \rightarrow \infty}
E \; =
    \di\sum_{s=1}^{k-1}
\lim_{n \rightarrow \infty}  \P(V_{s} < \eps) \P(
    s(\m n) = s)
& = 0.
\end{align*}
For the second quantity $F$ we obtain
\begin{align*}
\lim_{\m n \rightarrow \infty}   F
\  \leq \
      \di\sum_{s=k}^{v(K)}
  \lim_{\m n \rightarrow \infty}
\P(V_{s} < \eps )
\
 \leq \
(v(K)-k)
  \lim_{\m n \rightarrow \infty}
\P(V_{k} < \eps ),
\end{align*}
which is due to the fact that the statistics are embedded.
Let $(\ell,m)$ be such that
$r_{\cal V}(\ell,m)=k$. Since
$V_{k} > V_{D(\m n)}^{(\ell,m)}$, we have
\begin{eqnarray*}
  \lim_{\m n \rightarrow \infty}
\P(V_{k} < \eps )
& \leq &
  \lim_{\m n \rightarrow \infty}
\P(V_{D(\m n)}^{(\ell,m)} < \eps ).
 \end{eqnarray*}
Under $H_1(k)$, as in the proof of Theorem \ref{thm1_ind} we can see that the probability  $\P(D(\m n) < k)$ tends to zero as $\m n $ tends to infinity.
It follows that
$$ \lim_{\m n \rightarrow \infty}
\P(V_{D(\m n)}^{(\ell,m)} < \eps )
=\lim_{\m n \rightarrow \infty}
\P(V_{D(\m n)}^{(\ell,m)} < \eps \cap D(\m n) \geq  k)
$$
and since the statistics are embedded we have
$V_{k'}^{(\ell,m)} \geq \di \frac{n_{\ell}n_m}{n_{\ell}+n_m} \left(r_{\m j_0}^{(\ell,m)}\right)^2$ for all $k'\geq k$
which implies that
\beg
\label{limV}
\lim_{n \rightarrow \infty}
\P(V_{D(\m n)}^{(\ell,m)} < \eps )
&\leq&   \lim_{\m n \rightarrow \infty}
\P(\di \frac{n_{\ell}n_m}{n_{\ell}+n_m} \left(r_{\m j_0}^{(\ell,m)}\right)^2 < \eps )
\ \ = \ \ 0,
\en
since by (\ref{egalR00}) $  \di \frac{n_{\ell}n_m}{n_{\ell}+n_m} 
(r_{\m j_0 }^{(\ell,m)})^2
     =
    O_{\P}(\m n)$,
and finally
\begin{align*}
\lim_{n \rightarrow \infty}    \P(V_{s(\m n)} < \eps)
    \; \leq
    \lim_{n \rightarrow \infty} (E+F)
    &= 0.
        \end{align*}

\cqfd

%

\section{Paired  case}\label{appendIndep}
We briefly describe the adaptation in the case of dependent samples, rewriting the previous definitions and the main results. 
In what follows we write  $n:=n_1=\cdots = n_K$. 
%
\subsection{Two-sample paired  case}
The constructions
(\ref{Tkindep}) and (\ref{Tk2indep})
become
\begin{eqnarray*}\label{Tk_ind}
{T}^{(1, 2)}_{2,k}&=& \di n  \di\sum_{\m j\in{\cal S}(2); ord_{{}}(\m j,2)\leq k}   (r^{(1, 2)}_{\m j})^2,
{\rm \ for \ }1\leq k \leq  c(2),
\end{eqnarray*}
and, for $d>2$ and $1\leq k \leq c(d)$,
\begin{eqnarray*}\label{Tk2_ind}
{T}^{(1, 2)}_{d,k}&=& T^{(1, 2)}_{d-1,c(d-1)} + \di n  \di\sum_{\m j\in {\cal S}(d); ord_{{}}(\m j,d) \leq k }  ( r^{(1, 2)}_{\m j})^2.
\end{eqnarray*}
Then
(\ref{sumV}) and (\ref{rule})
become
\begin{align*}
    \label{sumV_ind}
V^{(1, 2)}_k &=
 \di n  \di \sum_{\m j \in{\cal  H}(k)}
  (r^{(1, 2)}_{\m j})^2
  \\
D(n)& := \min\big \{\argmax_{1\leq k \leq d(n)} (  V^{(1, 2)}_{k} - k q_{ n})  \big\},
\end{align*}
where
$q_{n}$  and $d(n)$ tend to $ +\infty$ as $n \to +\infty$.
A classical choice for $q_{ n}$ is $\alpha \log(n)$, where $\alpha$  can be simply equal to 1, or obtained by the tuning procedure described in Appendix \ref{annexTuning}.  

Finally, the associated  data-driven test statistic to  compare $C_{1}$ and $C_2$ is
\be*
\label{statpair_ind}
V^{(1, 2)}= V^{(1, 2)}_{D(n)}.
\e*

We consider the following rate for the penalty:

\begin{description}
\item{{\bf (A'')}}
$d(n)^{(p+5)} = o(p_{ n})$.
\end{description}
\begin{theorem}
\label{thm1}
If  {\bf (A'')} holds, then, under $ H_0$, ${D(n)}$ converges in Probability towards 1 as $n \rightarrow +\infty$.
\end{theorem}
Asymptotically, the null distribution will reduce to that of $V_{1}^{(1, 2)}=T^{(1, 2)}_{2,1}= n(r_{\m j}^{(1,2)})^2$, with $\m j = (1,1,0,\ldots,0)$ and 
\begin{align*}
    r^{(1, 2)}_{\m j}
    & =
    \di\frac{1}{n}
    \di\sum_{i=1}^{n}
\big(    
L_1(\w U^{(1)}_{i,1}) L_1(\w U^{(1)}_{i,2})
    -
   L_1(\w U^{(2)}_{i,1}) L_1(\w U^{(2)}_{i,2})
\big).
\end{align*}
\begin{theorem}\label{thm2_indE}
Let $\textbf{j}=(1,1,0\ldots,0)$. 
Then Under $H_0$,
\be*
\di (V^{(1, 2)})^{1/2}\overset{D}{\longrightarrow} \mathcal{N}\left(0,\sigma^2(1,2)\right),
\e*
where 

\be*
\di\sigma^2(1, 2)  &=& \V
\Bigg(
L_1(U^{(1)}_{1}) L_1(U^{(1)}_{2})
      -L_1(U^{(2)}_{1}) L_1( U^{(2)}_{2})
\\
&&
+ 2\sqrt{3}\di\int\di\int \big(\I(X^{(1)}_{1}\leq x )-F_{1}^{(1)}(x)\big)L_{1}( F^{(1)}_{2}(y))
d F^{(1)}(x,y)
\\
&&
- 2\sqrt{3}\di\int\di\int \big(\I(X^{(2)}_{1}\leq x )-F_{1}^{(2)}(x)\big)L_{1}( F^{(2)}_{2}(y))
d F^{(2)}(x,y)
\\
&&
+2\sqrt{3}\di\int\di\int \big(\I(X^{(1)}_{2}\leq y )-F_{2}^{(1)}(y)\big) L_{1}( F^{(1)}_{1}(x))
d F^{(1)}(x,y)
\\
&&
-
2\sqrt{3}\di\int\di\int \big(\I(X^{(2)}_{2}\leq y )-F_{2}^{(2)}(y)\big) L_{1}( F^{(2)}_{1}(x))
d F^{(2)}(x,y)
\Bigg).
\e*
\end{theorem}
To normalize the test, we consider the following estimator
\begin{align*}
\w \sigma^2{(1, 2)} & =
\di\frac{1}{n}
\di\sum_{i=1}^{n}\Big(M_{i,1} - M_{i,2} -\overline{M}_{1}+\overline{M}_{2}\Big)^2,
\end{align*}
 \begin{align*}
 \overline{M}_{s}  = \di\frac{1}{n}\sum_{i=1}^{n}M_{i,s},
 \ \
 \text{for \ } s=1,2,
\end{align*}
where
\begin{align*}
M_{i,s} &= L_{1}(\w U^{(s)}_{i,1})L_{1}(\w U^{(s)}_{i,2})
+
\di\frac{2\sqrt{3}}{n} \di\sum_{k=1}^{n}\Big( \I\big(X^{(s)}_{i,1}\leq X^{(s)}_{k,1}\big)-\w U^{(s)}_{k,1}\Big) L_{1}(\w U^{(s)}_{k,2})
\\
& \ \ \ 
+
\di\frac{2\sqrt{3}}{n} \di\sum_{k=1}^{n}\Big( \I\big(X^{(s)}_{i,2}\leq X^{(s)}_{k,2} \big)-\w U^{(s)}_{k,2}\Big) L_{1}(\w U^{(s)}_{k,1}).
\end{align*}
\begin{proposition}\label{prop1}
 Under $ H_0$,
\begin{align*}
\w\sigma^2(1, 2):= & \overset{\P}{\longrightarrow}
\sigma^2(1, 2).
\end{align*}
\end{proposition}
We then obtain the following result.
\begin{corollary} \label{cor1_ind}
Assume that  {\bf (A'')} holds. Under $H_0$,
$V^{(1, 2)}/\w \sigma^2(1, 2)$ converges in law towards  a chi-squared distribution $\chi_1^2$  as $n  \rightarrow +\infty$.
\end{corollary}

\subsection{K-sample paired  case}

The rule
(\ref{rule3})
becomes
\be* \label{rule3_ind}
s( n) & = & \min \Big\{\argmax_{1 \leq k \leq v(K)} \big({V}_{k} - k  p_{n}\big)  \Big\},
\e*
where $p_{ n}$ satisfies
\begin{description}
\item{{\bf (A''')}}
$d( n)^{p+5} = o(p_{ n})$.
\end{description}
In practice we choose $p_{n} =
\alpha \log(n)$.
We have 
\begin{theorem}
\label{thm2}
Assume that {\bf (A'')-(A''')} hold.  Under $ H_0$,
$s(n)$ converges in probability towards 1 as $n \rightarrow +\infty$.
\end{theorem}
\begin{corollary}\label{cor2_ind}
Assume that {\bf (A'')-(A''')} hold.
Under $ H_0$,
$V_{s(n)}/\w \sigma^2{(1,2)})$ converges in law towards a $\chi^2_1$ distribution.
\end{corollary}
Then the final data-driven test statistic is given by
\be*
V & = &
V_{s(n)}/\w \sigma^2{(1,2)}.
\e*

 \subsection{Alternative hypotheses in the paired case}
 We need the following assumption:
 \begin{description}
\item{{\bf (B')}}
$p_{ n}= o( n)$.
\end{description}
\begin{proposition}\label{thm4}
Assume that {\bf (A'')-(A''')-(B')} hold. Under $H_1(k)$, $s( n)$ converges in probability towards $k$ as $ n \rightarrow +\infty$, and
$V$ 
converges to $+\infty$, that is, $\P(V <\epsilon)\rightarrow 0$ for all $\epsilon >0$.
\end{proposition}
\begin{proposition}
\label{prop10}
Assume that {\bf (A'')} and  {\bf (A''')}  hold.  Then under $H_1^*$,
$s( n)$ converges in probability towards 1 as $n \rightarrow +\infty$ and
$V_{s( n)}/\w \sigma^2{(1,2)})$ converges in law towards a $\chi^2_1$ distribution.
\end{proposition}

\begin{proposition}\label{prop11}
Assume that {\bf (A'')-(A''')-(B')} hold. Then under $H_1^{**}(r_1,r_2)$, $s( n)$ converges in probability towards $r_1$ and $d( n)$ converges in probability towards $r_2$, as $n \rightarrow +\infty$. Moreover
$V$ 
converges to $+\infty$.
\end{proposition}
%
%

\section{Tuning the test statistic}
\label{annexTuning}

As evoked in Remark \ref{remark3}, we can choose the penalty \\ $q_{\m n}=p_{\m n} = \alpha \log(K^{(K-1)}n_1 \cdots n_K/(n_1+\cdots+n_K)^{K-1})$ by using the following   data-driven procedure.

~\\{\it Data-driven tuning procedure: 
\bit
\item
Assume we observe  $K$ populations, namely  $P_1, \ldots, P_K$
\item
Split randomly each population into $K'>2$ sub-populations, say
$P_{i,j}$, for $i=1,\ldots, K$, $j=1,\ldots, K'$.
\item
Clearly, for $i=1,\ldots, K$,   the $K'$ sub-populations $P_{i,1},\ldots, P_{i,K'}$ have the same copula, that is,  the null hypothesis $ H_0$ is satisfied.
\item
We can repeat $N$ times such a procedure to get $K*N$ $K'$ samples under the null.
\item
We then approximate numerically the value of the factor $\alpha>0$ such that the selection rule retains the first component, that is $s(\m n)=1$,  for all the   $K'$-sample tests.
From Theorem \ref{thm3_ind} this is the asymptotic expected value under the null. %

More precisely we fix
\begin{align*}
\w \alpha & =
\min \{ \alpha>0; \mbox{ such that  } s(\m n) =1   \mbox{ for the } K*N \mbox{selection rules} \}
\end{align*}
\eit
}

In our simulation, we fixed arbitrarily $K'=3$, which seems to give a very correct empirical level. Note that the use of this factor $\alpha$ only slightly modified the empirical results.

\section{Legendre polynomials}
\label{appendlegendre}
The Legendre polynomials used in this paper are defined on $[0,1]$ by
\be*
\label{legendre}
&&L_0=1,  L_1(x) = \sqrt{3}(2x-1) , {\rm \ and \ for \ } n>1:
\nonumber
\\
&&(n+1) L_{n+1}(x) =
\sqrt{(2n+1)(2n+3)}(2x-1)L_n(x)-\di\frac{n\sqrt{2n+3}}{\sqrt{2n-1}}L_{n-1}(x).
\e*
They satisfy
\be*
\di \int_0^1 L_j(x) L_k(x) dx & = & \delta_{jk},
\e*
where $\delta_{jk}=1$ if $j=k$ and $0$ otherwise.


\section{Representations of sepals and petals distributions}
\label{annexeIris}

\FloatBarrier

\begin{figure}[!ht]
	\begin{centering}
		 \includegraphics[height=7cm, width=7cm]{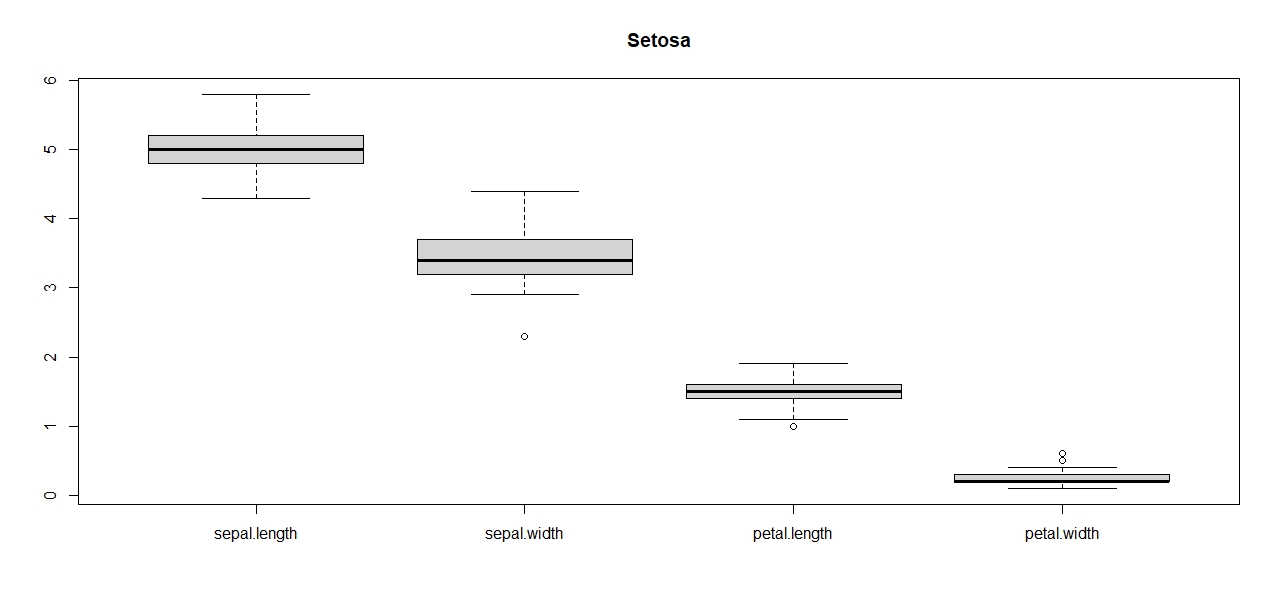}
		 \includegraphics[height=7cm, width=7cm]{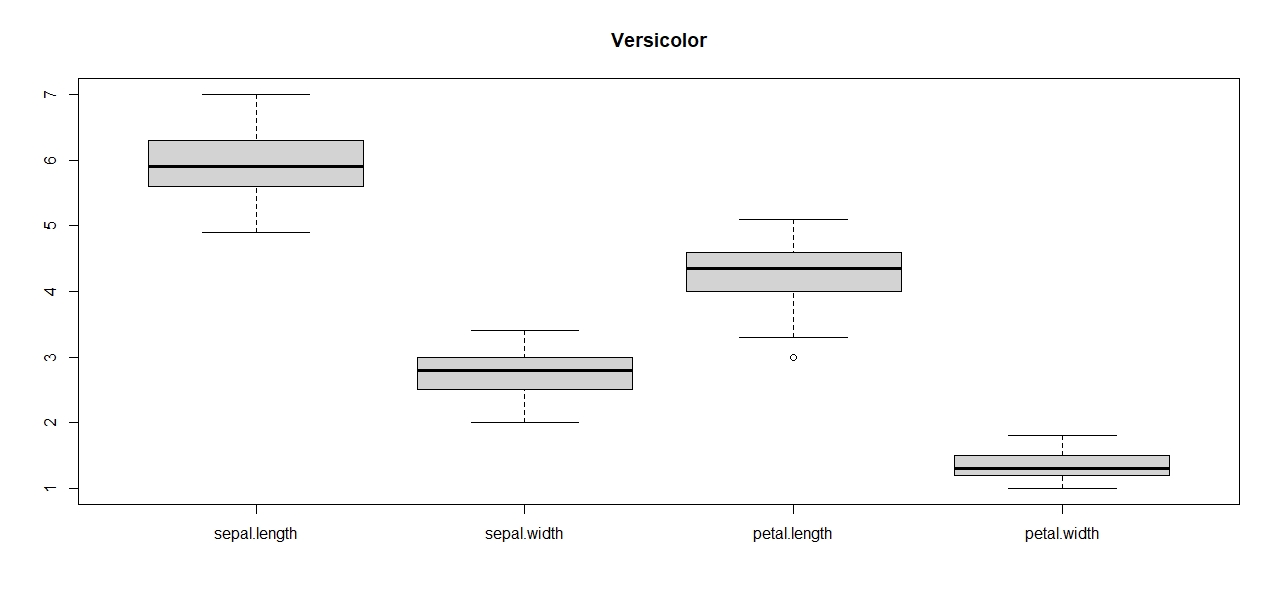}		
		 \includegraphics[height=7cm, width=7cm]{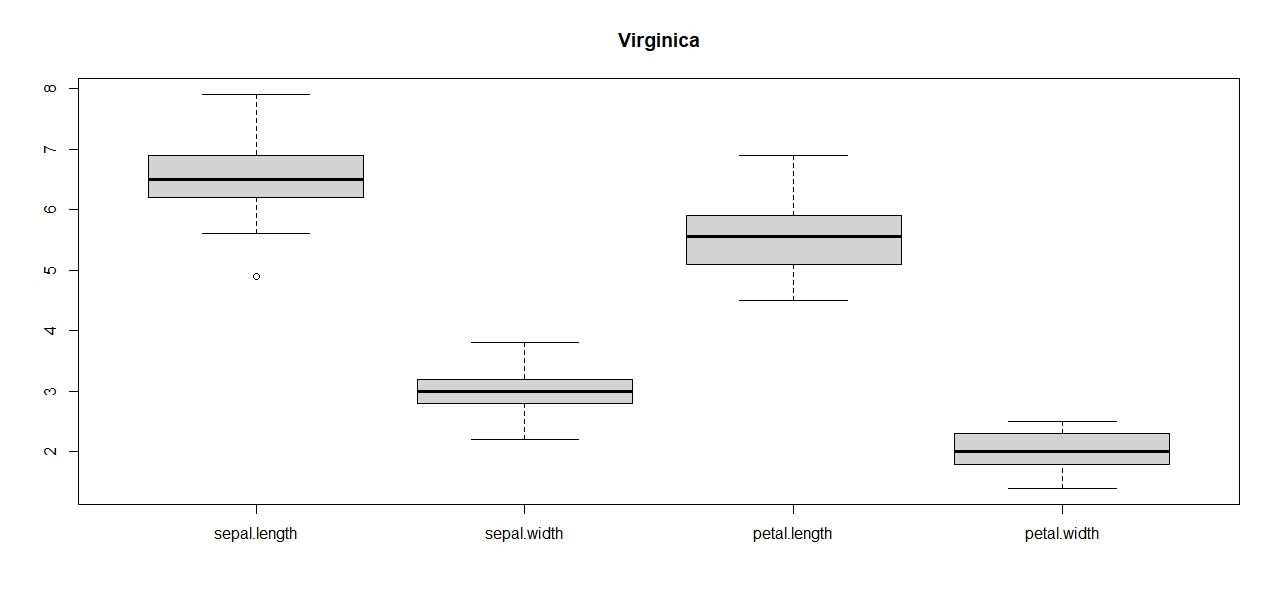}		
			\caption{Lengths and widths for Setosa,  Versicolor and Virginica.}
			\label{figiris}
	\end{centering}
\end{figure}

\FloatBarrier

\section{Simulation results in the two-sample case (complements)}
\label{appendRemillard}


\begin{figure}[!htbp]
    \centering 
\begin{subfigure}{0.49\textwidth}
  \includegraphics[width=\linewidth]{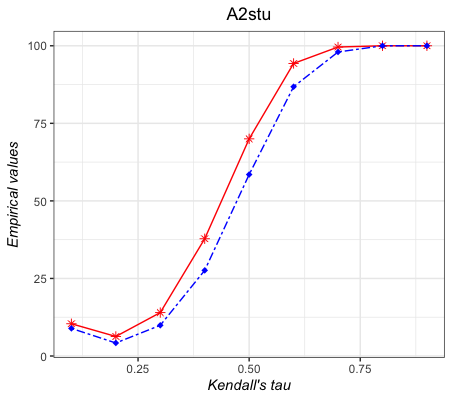}
\end{subfigure}\hfil 
\begin{subfigure}{0.49\textwidth}
  \includegraphics[width=\linewidth]{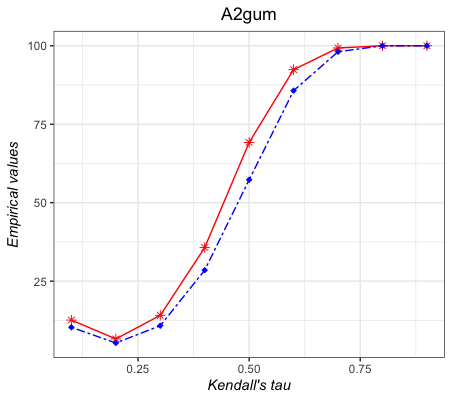}
\end{subfigure}

\medskip
\begin{subfigure}{0.49\textwidth}
  \includegraphics[width=\linewidth]{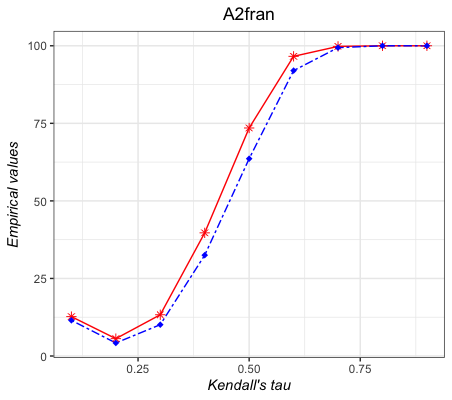}
\end{subfigure}\hfil 
\begin{subfigure}{0.49\textwidth}
  \includegraphics[width=\linewidth]{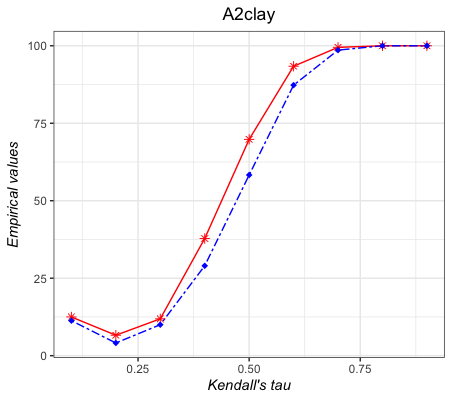}
\end{subfigure}\hfil 
\begin{subfigure}{0.49\textwidth}
  \includegraphics[width=\linewidth]{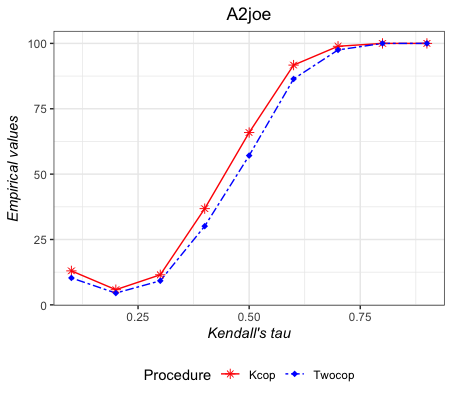}
\end{subfigure}
\caption{Two-sample case: empirical powers under  alternatives $\mathcal{A}2:\textbf{50-50}$.}
\label{powerA25050}
\end{figure}

\FloatBarrier
\begin{figure}[!htbp]
    \centering 
\begin{subfigure}{0.49\textwidth}
  \includegraphics[width=\linewidth]{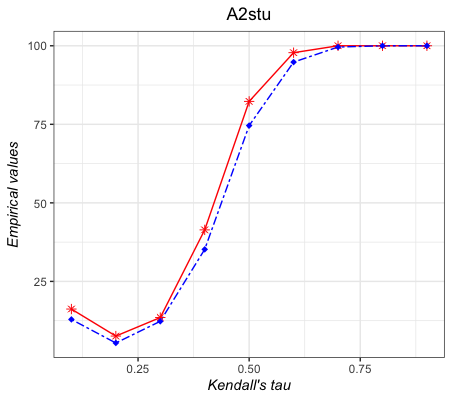}
\end{subfigure}\hfil 
\begin{subfigure}{0.49\textwidth}
  \includegraphics[width=\linewidth]{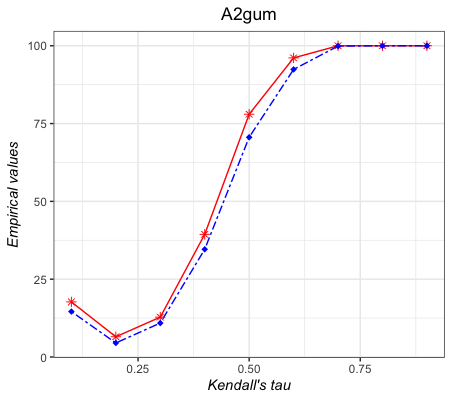}
\end{subfigure}

\medskip
\begin{subfigure}{0.49\textwidth}
  \includegraphics[width=\linewidth]{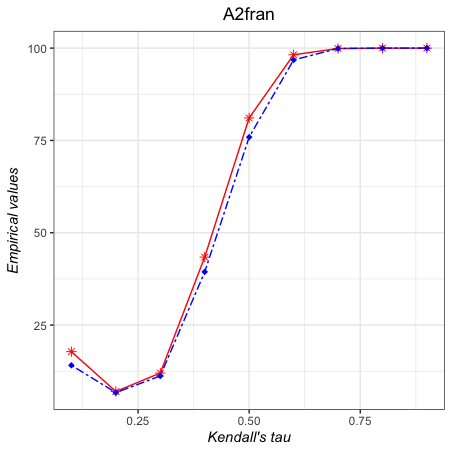}
\end{subfigure}\hfil 
\begin{subfigure}{0.49\textwidth}
  \includegraphics[width=\linewidth]{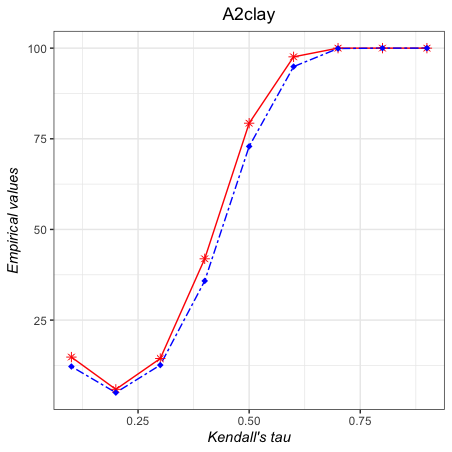}
\end{subfigure}\hfil 
\begin{subfigure}{0.49\textwidth}
  \includegraphics[width=\linewidth]{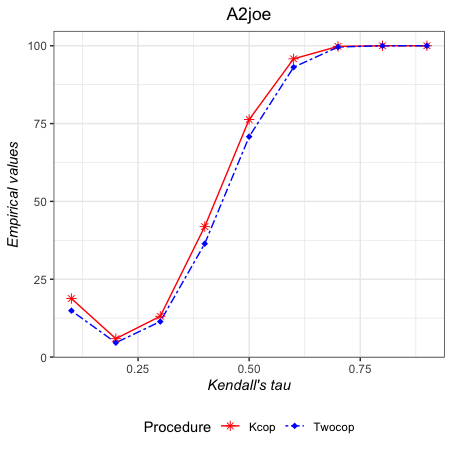}
\end{subfigure}
\caption{Two-sample case: empirical powers under  alternatives $\mathcal{A}2:\textbf{50-100}$.}
\label{powerA250100}
\end{figure}

\FloatBarrier
\begin{figure}[!htbp]
    \centering 
\begin{subfigure}{0.49\textwidth}
  \includegraphics[width=\linewidth]{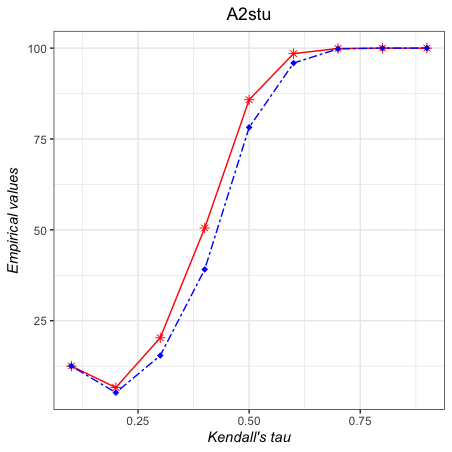}
\end{subfigure}\hfil 
\begin{subfigure}{0.49\textwidth}
  \includegraphics[width=\linewidth]{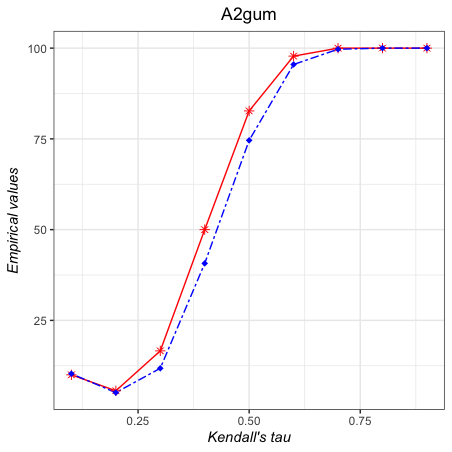}
\end{subfigure}

\medskip
\begin{subfigure}{0.49\textwidth}
  \includegraphics[width=\linewidth]{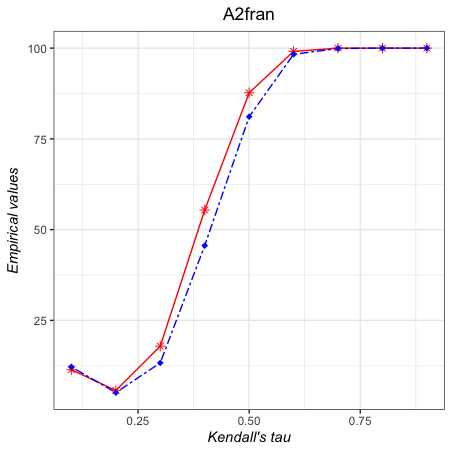}
\end{subfigure}\hfil 
\begin{subfigure}{0.49\textwidth}
  \includegraphics[width=\linewidth]{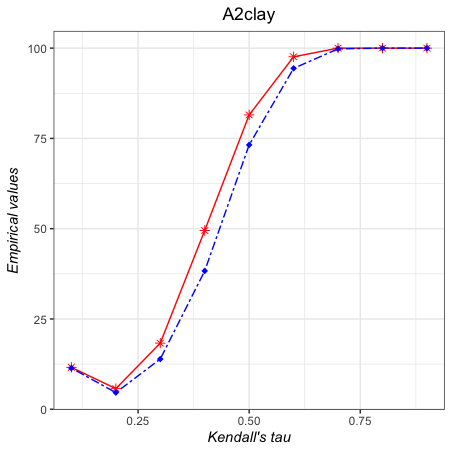}
\end{subfigure}\hfil 
\begin{subfigure}{0.49\textwidth}
  \includegraphics[width=\linewidth]{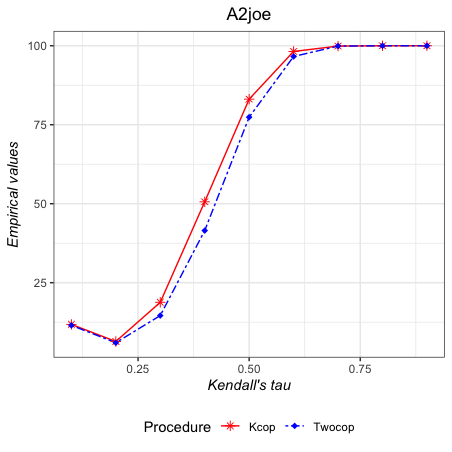}
\end{subfigure}
\caption{Two-sample case: empirical powers under  alternatives $\mathcal{A}2:\textbf{100-50}$.}
\label{powerA210050}
\end{figure}

\FloatBarrier
\begin{figure}[!htbp]
    \centering 
\begin{subfigure}{0.49\textwidth}
  \includegraphics[width=\linewidth]{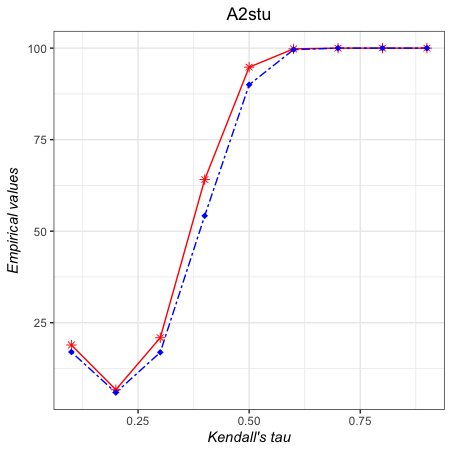}
\end{subfigure}\hfil 
\begin{subfigure}{0.49\textwidth}
  \includegraphics[width=\linewidth]{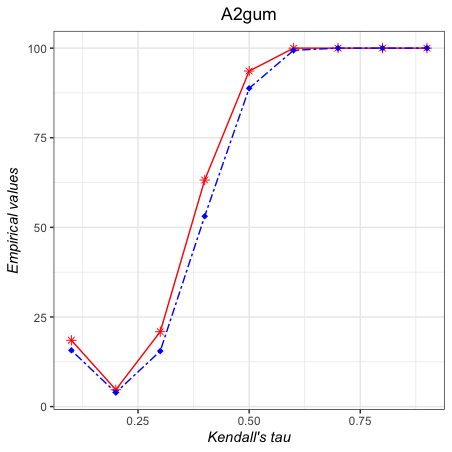}
\end{subfigure}

\medskip
\begin{subfigure}{0.49\textwidth}
  \includegraphics[width=\linewidth]{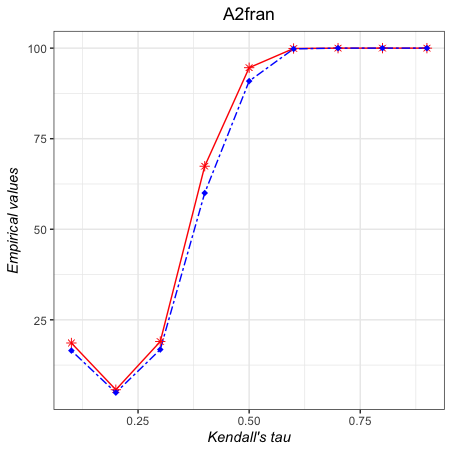}
\end{subfigure}\hfil 
\begin{subfigure}{0.49\textwidth}
  \includegraphics[width=\linewidth]{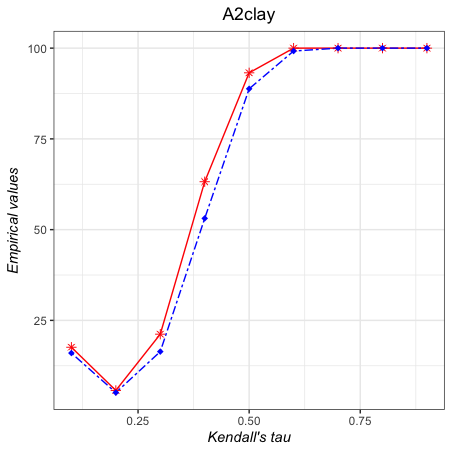}
\end{subfigure}\hfil 
\begin{subfigure}{0.49\textwidth}
  \includegraphics[width=\linewidth]{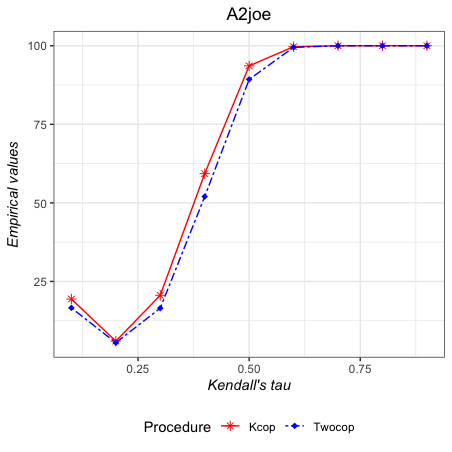}
\end{subfigure}
\caption{Two-sample case: empirical powers under alternatives  $\mathcal{A}2:\textbf{100-100}$.}
\label{powerA2100100}
\end{figure}

\FloatBarrier

\section{Insurance data: the two-by-two comparison}
\label{appendANOVA}
\FloatBarrier
%
%


\begin{table}[htbp!]
\caption{ANOVA test p-values (in bold the cases where the equality is not rejected). Values given in brackets indicate the size of groups ($G$).}
\label{anova_pvalue}
\begin{tabular}{@{}llllllllll@{}}
\hline
Groups  & $G_1$ & $G_2$ & $G_3$ & $G_4$ & $G_5$ & $G_6$ & $G_7$ & $G_8$ & $G_9$ \\\hline
$G_1(966)$ & & & & & & & & & \\
$G_2(971)$ & {\bf 0.794} & & & & & & & & \\
$G_3(996)$ & {\bf 0.265} &{\bf  0.193} & & & & & & & \\
$G_4(954)$ & {\bf  0.827} & {\bf 0.952} & {\bf  0.175} & & & & & & \\
$G_5(955)$ & {\bf  0.397} & {\bf  0.588} & {\bf 0.051} & {\bf  0.519} & & & & & \\
$G_6(915)$ & {\bf 0.066} & {\bf 0.138} & 0.003 & {\bf  0.10} & {\bf  0.325} & & & & \\
$G_7(828)$ & 0.002 & 0.009 & 0.000 & 0.005 & 0.028 & {\bf  0.209} & & & \\
$G_8(624)$ & 0.001 & 0.005 & 0.000 & 0.002 & 0.017 & {\bf  0.152} &{\bf  0.883} & & \\
$G_9(524)$ & 0.030 & {\bf 0.069} & 0.001 & 0.046 & {\bf  0.179} & {\bf  0.700} & {\bf  0.389} & {\bf 0.304} & \\
$G_{10}(396)$ &  0.008 & 0.020 & 0.000 & 0.013 & {\bf  0.056} & {\bf  0.289} & {\bf  0.925} & {\bf  0.816} &{\bf  0.483} \\\hline
\end{tabular}
\end{table}

\section{3-sample with Student copulas}
\label{appendstudent}
We consider three Student copulas  
$C_1,C_2, C_3$,   with df=5 and Kendall's tau $\tau_1,\tau_2,\tau_3$, respectively. The first alternative is a very smooth deviation   $(\tau_1,\tau_2,\tau_3)$ = $(0.4, 0.5, 0.6)$ coinciding with three closely related populations. The second alternative is formed of only two populations but with a slightly larger difference $(\tau_1,\tau_2,\tau_3)= (0.4, 0.4, 0.6)$.  
Table \ref{6poptab} contains empirical powers for $n\in\{50,100, 200\}$. 
%
It appears to be easier to detect the second alternative, which involves two more distinct groups, rather than three groups with a smooth variation. This finding may suggest the possibility of a forward test-based clustering procedure, wherein each population is successively tested before being joined to a cluster. This perspective could be explored further.
\FloatBarrier

\begin{table}[htbp!]
\caption{Empirical powers}
\label{6poptab}
\begin{tabular}{cccc}
\hline
  & $n=50$  & $n=100$ & $n=200$  \\
$(\tau_1,\tau_2,\tau_3)=(0.4,0.4,0.6)$ & 56.4& 73.2 & 96.6\\
$(\tau_1,\tau_2,\tau_3)=(0.4,0.5,0.6)$ & 17.6 & 29.8 & 39.5\\
\end{tabular}
\end{table}

%
%
%
%
%
%



\section{Empirical levels for the ten-sample case}

\label{appendempirical}

\begin{table*}[ht]
\caption{Empirical levels for the ten-sample test.}
\centering
\begin{tabular}{ccccccc}
\hline
& \multicolumn{6}{c}{Models} \\ \cline{2-7}
$n$  & Gaussian & Student  & Gumbel & Frank  & Clayton & Joe \\
  \hline 
& \multicolumn{6}{c}{Kendall tau $\tau=0.1$} \\ \cline{2-7}
50 & 11.5 & 12.1 & 10.6 & 10.9 & 11.0 & 10.8 \\
100 & 9.9 & 9.3 & 9.3 & 9.6 & 8.3 & 8.3 \\
200 & 7.8 & 6.2 & 7.9 & 6.2 & 7.5 & 7.8 \\
300 & 6.9 & 7.5 & 7.0 & 7.0 & 5.7 & 6.8 \\
400 & 6.4 & 5.1 & 6.7 & 5.7 & 5.3 & 6.0 \\
500 & 5.2 & 6.0 & 5.9 & 6.2 & 7.1 & 5.7 \\
600 & 5.6 & 7.4 & 5.2 & 6.4 & 5.7 & 5.6 \\
700 & 5.1 & 6.3 & 5.4 & 6.0 & 5.2 & 7.2 \\
800 & 5.1 & 5.6 & 6.2 & 5.8 & 6.3 & 5.8 \\
900 & 5.8 & 3.4 & 5.6 & 6.2 & 5.3 & 6.6 \\
1000 & 6.0 & 5.9 & 5.1 & 4.2 & 6.4 & 5.1 \\
\cline{2-7}
& \multicolumn{6}{c}{Kendall tau $\tau=0.5$} \\ \cline{2-7}
 50 & 5.4 & 4.0 & 3.6 & 3.2 & 4.4 & 3.7 \\
100 & 6.0 & 4.2 & 4.0 & 5.4 & 5.6 & 3.6 \\
200 & 4.9 & 4.5 & 5.2 & 5.1 & 5.3 & 4.3 \\
300 & 5.7 & 5.5 & 5.5 & 4.7 & 5.4 &4.0 \\
400 & 4.7 & 5.0 & 5.4 & 4.6 & 5.3 & 3.2 \\
500 & 4.4 & 4.9 & 4.1 & 5.5 & 5.5 & 4.8 \\
600 & 4.8 & 6.5 & 5.1 & 6.1 & 4.8 & 6.2 \\
700 & 5.4 & 5.2 & 6.1 & 4.6 & 4.8 & 3.9 \\
800 & 4.9 & 6.3 & 4.5 & 6.1 & 4.9 & 4.8 \\
900 & 4.6 & 4.0 & 4.8 & 5.2 & 4.9 & 4.2 \\
1000 & 4.2 & 5.5 & 4.5 & 4.1 & 4.8 & 3.6 \\

\cline{2-7}
& \multicolumn{6}{c}{Kendall tau $\tau=0.8$} \\ \cline{2-7}

50 & 1.0 & 0.6 & 0.6 & 0.8 & 3.1 & 0.0 \\
100 & 2.6 & 2.5 & 1.8 & 2.0 & 4.8 & 0.7 \\
200 & 4.1 & 4.0 & 4.3 & 4.0 & 5.1 & 2.3 \\
300 & 4.0 & 4.5 & 3.6 & 4.0 & 5.7 & 4.3 \\
400 & 3.5 & 4.1 & 4.9 & 3.7 & 5.0 & 3.3 \\
500 & 4.9 & 3.9 & 3.6 & 4.8 & 3.9 & 4.4 \\
600 & 4.6 & 5.2 & 5.7 & 4.9 & 4.9 & 4.8 \\
700 & 4.0 & 5.0 & 5.5 & 4.9 & 4.6 & 4.0 \\
800 & 4.5 & 6.5 & 3.2 & 3.7 & 4.3 & 3.5 \\
900 & 4.4 & 4.6 & 4.0  & 5.9 & 5.8 & 4.5 \\
1000 & 3.7 & 5.5 & 4.7 & 4.3 & 5.2 & 4.7 \\

\hline
\end{tabular}
\end{table*}


\end{document}